\documentclass[12pt]{article}
\usepackage{amsmath,amsfonts,amssymb,amsthm,mathdots}
\usepackage{authblk}
\usepackage{color}
\usepackage{comment}
\usepackage[normalem]{ulem}
\usepackage{fullpage}
\usepackage{xcolor}
\usepackage{hyperref}

\numberwithin{equation}{section}

\def\Z{\mathbb{Z}}

\def\F{\mathbb{F}}

\newcommand{\pres}[2]{\langle {#1}\ |\ {#2} \rangle}

\DeclareMathOperator{\cont}{cont}
\DeclareMathOperator{\res}{Res}
\DeclareMathOperator{\adj}{adj}

\newtheorem{theorem}{Theorem}[section]
\newtheorem{corollary}[theorem]{Corollary}
\newtheorem{proposition}[theorem]{Proposition}
\newtheorem{lemma}[theorem]{Lemma}
\newtheorem{example}[theorem]{Example}
\theoremstyle{definition}

\newtheorem*{claim*}{Claim}

\theoremstyle{definition}
\newtheorem{remark}[theorem]{Remark}

\title{Smith forms of matrices in Companion Rings, \linebreak with group theoretic and topological applications}

\author{Vanni Noferini and Gerald Williams}

\begin{document}
    
\newcommand{\Addresses}{{
  \bigskip
  \footnotesize

  Vanni Noferini (corresponding author), \textsc{Department of Mathematics and Systems Analysis, Aalto University PL 11000, 00076 Aalto, Finland}\par\nopagebreak
  \textit{E-mail address}: \texttt{vanni.noferini@aalto.fi}

  \medskip

  Gerald Williams, \textsc{School of Mathematics, Statistics and Actuarial Science, University of Essex, Colchester, Essex CO4 3SQ, United Kingdom}\par\nopagebreak
  \textit{E-mail address}: \texttt{gerald.williams@essex.ac.uk}

}}

\maketitle

\begin{abstract}
Let $R$ be a commutative ring and $g(t) \in R[t]$ a monic polynomial.  The commutative ring of polynomials $f(C_g)$ in the companion matrix $C_g$ of $g(t)$, where $f(t)\in R[t]$, is called the Companion Ring of $g(t)$. Special instances include the rings of circulant matrices, skew-circulant matrices, pseudo-circulant matrices, or lower triangular Toeplitz matrices.
When $R$ is an Elementary Divisor Domain, we develop new tools for computing the Smith forms of matrices in Companion Rings. In particular, we obtain a formula for the second last non-zero determinantal divisor, we provide an $f(C_g) \leftrightarrow g(C_f)$ swap theorem, and a composition theorem. When $R$ is a principal ideal domain we also obtain a formula for the number of non-unit invariant factors. By applying these to families of circulant matrices that arise as relation matrices of cyclically presented groups, in many cases we compute the groups' abelianizations. When the group is the fundamental group of a three dimensional manifold, this provides the homology of the manifold.  In other cases we obtain lower bounds for the rank of the abelianisation and record consequences for finiteness or solvability of the group, or for the Heegaard genus of a corresponding manifold.
\end{abstract}

\noindent {\bf  MSC:} 15A21, 15B36, 11C99, 20F05, 57M05, 57M10.

\noindent {\bf Keywords:} Smith form; Elementary Divisor Domain; companion matrix; Companion Ring; circulant; cyclically presented group; Fibonacci group; abelianization; homology.

\section{Introduction}

Companion Rings were introduced and studied in \cite{NW21}. These are rings of matrices $f(C_g)$, whose entries are elements of a commutative ring $R$, and where $f(t),g(t)$ are polynomials in $R[t]$, where $g(t)$ is monic with companion matrix $C_g$. They include, as special cases, the rings of circulant, skew-circulant, pseudo-circulant, and lower triangular Toeplitz matrices. (Formal definitions will be given in Section \ref{sec:background}.) The Smith form of a matrix over an Elementary Divisor Domain  (EDD) $R$ \cite{Friedland} is a diagonal matrix whose diagonal entries are elements of $R$ that form a divisor chain; the Smith form of a matrix can be expressed in terms either of the invariant factors or of the determinantal divisors of the matrix. In \cite{NW21} properties of Smith forms of matrices in Companion Rings were obtained. This includes formulae for the number of non-zero determinantal divisors, for the first determinantal divisor, and the last non-zero determinantal divisor, as well as theorems concerning factorizations of $f(t)$ and $g(t)$. The results obtained were designed for the primary application sought in that article, namely, the computation of the first integral homology of 3-dimensional Brieskorn manifolds.

In this article we build on \cite{NW21} by developing further tools, with a wider range of applications, for calculating the Smith form of a matrix $f(C_g)$ in a Companion Ring. In particular, we obtain formulae for the second last non-zero determinantal divisor (Theorem \ref{thm:gamman-1}), a ``swap theorem'' that allows us to interchange between $f(C_g)$ and $g(C_f)$ for monic polynomials $f(t)$ and $g(t)$ (Theorem \ref{thm:swaptheorem}), and a composition theorem concerning matrices of the form $(f\circ h)(C_{g\circ h})$ where $h(t)$ is a monic polynomial in $R[t]$ (Theorem \ref{thm:composition}). In the case when $R$ is a Principal Ideal Domain (PID), we obtain a formula for the number of non-unit invariant factors of $f(C_g)$ (Theorem \ref{thm:genJO}), recovering, as a particular instance, a result of Johnson and Odoni \cite{JohnsonOdoni}.

We apply our results to families of circulant matrices. While, in principle, they could be applied to any situation in which circulant matrices play a role, our applications focus on circulant matrices that arise as relation matrices of various classes of cyclically presented groups. These are groups defined in terms of a group presentation with an equal number of generators and relators that admits a cyclic symmetry. The relation matrix -- a matrix whose  Smith form provides the abelianization of the group -- is a circulant matrix. These applications are in a similar spirit to \cite{NW22} which used the results of \cite{NW21} to identify classes of cyclically presented groups with free abelianization. Here, we go further, by calculating completely the abelianization (which in certain cases provides the homology of a corresponding 3-dimensional manifold), or by obtaining a lower bound for the number of generators of the abelianization. In turn, this has consequences concerning the finiteness or solvability of the group, or provides an application to low-dimensional topology by delivering a lower bound for the Heegaard genus of the manifold.

Specifically, our applications are as follows. In Section \ref{sec:circulantgraphs} we give a new, short, proof of the (already known) Smith form of the adjacency matrix of the Cocktail party graphs. In Section \ref{sec:FractionalFibonacci} we obtain explicit formulae for the abelianization of the fractional Fibonacci groups; this result has been previously stated (without proof) for particular choices of the defining parameters, but the general case is new.  In Section \ref{sec:neuwirth} we obtain the Smith form of the circulant matrix whose rows are cyclic permutations of the vector $[a,\dots,a, b, \dots, b]$, where the number of $b$'s is coprime to the order of the matrix; in particular, we recover the special case where each row contains exactly one $b$  which was stated (without proof) in \cite{SpaggiariNeuwirth} and provides the homology of the periodic generalized Neuwirth manifolds.
The remaining applications are all completely new results. In Section \ref{sec:Hrns} we obtain formulae for the abelianizations of generalised Fibonacci group $\mathcal{H}(r,n,s)$ under certain conditions on the parameters. In Section \ref{sec:positivelength3} we consider certain cyclically presented groups whose relators are positive words of length three, and we prove two results conjectured in \cite{MohamedWilliams}: we show that the abelianizations of two classes of such groups are isomorphic, and we compute the rank of the abelianization of groups in one such family. In Section \ref{sec:CRSgroups} we obtain a sharp lower bound for the minimum number of generators of an $8$-parameter class of cyclically presented group that encompasses many classes of groups that arise in topological settings. This result has implications for the finiteness and solvability of the groups, and for the Heegaard genus of a corresponding manifold. 

\section{Preliminaries}\label{sec:background}

\subsection{The Smith Theorem}

Given a commutative ring $R$, a \emph{unimodular} matrix \cite[p.\,12]{Newmanbook} is a square matrix $U \in R^{n \times n}$ such that $\det U$ is a unit of $R$. Unimodular matrices are precisely the units of $R^{n \times n}$, that is, matrices whose inverse exists in $R^{n \times n}$. The Smith Theorem, proved by Smith in \cite{Smith} for the case $R=\Z$, and later by Kaplansky \cite{kap49} for EDDs, is stated as Theorem \ref{thm:SNF} below. See, for example, \cite[Theorem 1.14.1]{Friedland} for a proof.

\begin{theorem}[Smith Theorem]\label{thm:SNF}
Let $R$ be an EDD and $M \in R^{m \times n}$. Then there exist unimodular matrices $U \in R^{n \times n}$, $V \in R^{m \times m}$ such that $U M V = S$ where $S$ is diagonal and satisfies $S_{i,i} \mid S_{i+1,i+1}$ for all $i=1,\dots,\min(m,n)-1$. Further, let $\gamma_0(M)=1 \in R$, and for $i=1,\dots,\min(m,n)$ define the $i$-th \emph{determinantal divisor} $\gamma_i(M)$ to be the greatest common divisor GCD of all minors of $M$ of order $i$. Then
$$S_{i,i} = \frac{\gamma_i(M)}{\gamma_{i-1}(M)} =: s_i(M),$$
where the diagonal elements $s_i(M)$, $i=1,\dots,\min(m,n)$, are called the \emph{invariant factors} of $M$. The matrix $S$ is called the \emph{Smith form} of $M$.
\end{theorem}

When the matrix $M$ in question is clear from the context we will simply write $\gamma_i$ or $s_i$ rather than $\gamma_i(M)$ or $s_i(M)$.

\begin{remark}
\begin{itemize}
    \item[(a)] The GCDs of subsets of a ring $R$ are defined up to multiplication by units of $R$, and so the invariant factors, determinantal divisors, and Smith form of a matrix are defined up to multiplication by units (although often a convention is set such as, in the case $R=\Z$, that GCDs are positive integers, to impose uniqueness). When we refer to these objects, we will implicitly mean that they are taken up to units of the ring. Similarly, when we refer to resultants, we use the symbol $\res(f(t), g(t))$ to mean any associate of the resultant of $f(t)$ and $g(t)$.

    \item[(b)] We stress that, in this paper, we find it convenient to define an invariant factor to be any, and possibly zero, diagonal element of the Smith form. This is in contrast with the, common in the literature (see e.g. \cite[p.\,28]{Newmanbook}), convention to call invariant factors only the \emph{nonzero} diagonal elements in the Smith form.
\end{itemize}
\end{remark}

Two matrices $M,N\in R^{m\times n}$ are said to be \emph{equivalent} (over $R$), denoted $M\sim N$, if there exist unimodular matrices $U \in R^{m \times m}, V \in R^{n \times n}$ such that $UMV=N$, and they are said to be \emph{similar} (over $R$), denoted $M\sim_S N$, if there exists a unimodular matrix $U$ such that $UMU^{-1}=N$.  It follows from the Smith Theorem that any pair of $m \times n$ matrices with entries in $R$ are equivalent if and only if they have the same invariant factors and that, since rank is preserved by multiplication by invertible matrices, $M$ has rank $r$ if and only if its invariant factors satisfy $s_i(M)=0$ precisely when $i>r$.

We conclude this subsection by recalling the definition and some properties of the adjugate (sometimes called adjoint) of a square matrix $M \in R^{n \times n}$ where $R$ is an integral domain; see \cite[Section 0.8.2]{HJ85} for more details. The \emph{adjugate} of $M$ is the matrix $\adj(M) \in R^{n \times n}$ whose $(i,j)$ entry is $(-1)^{i+j} \mathcal{M}_{ij}$ where $\mathcal{M}_{ij} \in R$ is the determinant of the $(n-1) \times (n-1)$ matrix obtained by removing the $i$-th row and the $j$-th column of $M$. It follows from the definition that $\gamma_{n-1}(M)=\gamma_1(\adj(M))$ is the GCD of the entries of $\adj(M)$; another useful property that we will freely use in the paper is $\adj(M) M = M \adj(M) = \det(M) I_n$. 

\subsection{Scalar polynomials in $R[t]$}

The following theorem (see for example \cite[pp.\,128--129]{Jacobson}) concerns division in the polynomial ring $R[t]$, where $R$ is a commutative ring with unity. We adopt the convention $\deg 0 := -\infty$. With this convention one has $\deg a(t)b(t) = \deg a(t) + \deg b(t)$ and $\deg (a(t) + b(t)) \leq \max \{ \deg a(t), \deg b(t) \}$ (with exact equality if the maximum is only attained once).

\begin{theorem}[Polynomial division]\label{thm:quotientremainder}
    Let $R$ be a commutative ring with unity and suppose that $g(t) \in R[t]$ is such that the leading coefficient of $g(t)$ is a unit of $R$. Then, polynomial division by $g(t)$ is well defined in $R[t]$, that is, for every $f(t) \in R[t]$ there exist unique polynomials $q(t),r(t)$ (called, resp., the \emph{quotient} and \emph{remainder} of the division of $f(t)$ by $g(t)$) such that $f(t)=g(t)q(t)+r(t)$ and $\deg r(t) < \deg g(t)$.
\end{theorem}

If $R$ is a GCD domain, the \emph{content} of a nonzero polynomial $f(t)=\sum_{i=0}^m f_i t^i \in R[t]$ is the GCD of the coefficients of $f(t)$, i.e., $\mathrm{cont}(f(t))=\gcd(f_0,\dots,f_m)$; the content of the zero polynomial is defined to be $\mathrm{cont}(0)=0$.

\subsection{Matrices in Companion Rings}\label{sec:companionrings}

Let $R$ be a commutative ring (with unity) other than $\{0\}$, and fix the monic polynomial $g(t)=t^n + \sum_{i=0}^{n-1} g_i t^i \in R[t]$. The ideal $\langle g(t) \rangle \subset R[t]$ is the set of polynomials that are multiples of $g(t)$. Moreover, we write $a(t) \equiv b(t) \bmod g(t)$ if $a(t)-b(t) \in \langle g(t) \rangle$; this notation extends elementwise to matrices, i.e., we write $A(t) \equiv B(t) \bmod g(t)$ if every entry of $A(t)-B(t)$ is divisible by $g(t)$. The quotient ring $\mathcal{Q}:=R[t]/\langle g(t) \rangle$ is the set of equivalence classes with respect to the above defined equivalence mod $g(t)$: namely, an element of $\mathcal{Q}$ has the form $[f(t)]=\{a(t) \in R[t] \ :\ a(t) \equiv f(t) \bmod g(t)\}$. Furthermore, by Theorem \ref{thm:quotientremainder} and since $g(t)$ has leading coefficient $1$, polynomial division by $g(t)$ uniquely defines a quotient and a remainder. In particular, this means that each equivalence class in $\mathcal{Q}$ has a unique representative having degree strictly less than $n$. It follows that $\mathcal{Q}$ is a module over $R$; a (canonical) basis of $\mathcal{Q}$ is the monomial basis $\{[1],[t],\dots,[t^{n-1}] \}$. Let us consider the linear endomorphism  
\[ M_{[t]} : \mathcal{Q} \rightarrow \mathcal{Q}, \quad [a(t)] \mapsto [t a(t)].\]
The representation of $M_{[t]}$ in the canonical basis is denoted by $C_g$ and it is called the \emph{companion matrix} of the polynomial $g(t)$. Explicitly,
\begin{equation}\label{eq:cg}
    C_g = \begin{bmatrix}
        -g_{n-1} & \dots & -g_1 & -g_0\\
        1 & & & \\
        & \ddots & & \\
        & & 1 &
    \end{bmatrix} \in R^{n \times n},
\end{equation}
where entries not explicitly displayed are understood to be $0$. For any integer $k\geq 1$ define
\begin{equation}
 \Lambda_k(t)=\begin{bmatrix}
    t^{k-1} & t^{k-2} & \dots & t & 1
\end{bmatrix}^T.\label{eq:Lambdan(t)}
\end{equation}
Then it can be readily verified that
\[ C_g \Lambda_n(t) = t \Lambda_n(t) - g(t) \begin{bmatrix}
1 & 0& \dots & 0
\end{bmatrix}^T, \]
thus showing that $C_g$ represents $M_{[t]}$ in the monomial basis of $\mathcal{Q}$. Indeed, $C_g$ is defined uniquely by the property $C_g \Lambda(t) \equiv t \Lambda(t) \bmod g(t)$ (taking into account that the monomial basis has been chosen to represent elements of $\mathcal{Q}$). 

Let us now consider the (commutative!) subring of $R^{n \times n}$ consisting of matrices of the form $f(C_g)=\sum_i f_i C_g^i$ where $f(t)=\sum_{i} f_i t^i \in R[t]$ is a polynomial with coefficients in $R$; following \cite{NW21}, we call this ring the \em Companion Ring \em of $g(t)$. It follows from the form of $C_g$ that for each $0\leq k<n$ the bottom row of $C_g^k$ is $e_{n-k}^T$. Hence, if $f(t)=\sum_{i=0}^{n-1} f_i t^i$ has degree strictly less than $n$, then the bottom row of $f(C_g)$ is $[f_{n-1},\dots ,f_1,f_0]$. In addition, the defining property $C_g \Lambda(t) \equiv t \Lambda(t) \bmod g(t)$ extends to 
 polynomial functions of $C_g$: indeed, it is a consequence of the third isomorphism theorem for rings that $\mathcal{Q}$, $R^{n}$, and the Companion Ring of $g(t)$ are all isomorphic. Generally, for any $f(t) \in R[t]$ the map $M_{[f(t)]}:[a(t)] \mapsto [a(t)f(t)]$ is represented (in the monomial basis) by the matrix $f(C_g)$. Thus, by the observations above, $f(C_g)$ is the unique element of $R^{n \times n}$ that satisfies $f(C_g) \Lambda(t) \equiv f(t) \Lambda(t) \bmod g(t)$.  We state this property formally below, as we use it frequently throughout the paper:

\begin{proposition}\label{prop:trick}
Let $R$ be a commutative ring with unit and such that $0 \neq 1$, let $g(t)=t^n + \sum_{i=0}^{n-1} g_i t^i \in R[t]$ be monic, $\Lambda_n(t)$ be as in \eqref{eq:Lambdan(t)} and let $C_g \in R^{n \times n}$ be the companion matrix of $g(t)$ as in \eqref{eq:cg}.
For any polynomial $f(t)\in R[t]$, if $X \in R^{n \times n}$ satisfies $X \Lambda_n(t) \equiv f(t) \Lambda_n(t) \bmod g(t)$, then $X=f(C_g)$.
\end{proposition}

\begin{remark}
    Note that, by the Cayley-Hamilton theorem and for all $f(t) \in R[t]$, $f(C_g)=\phi(C_g)$ where $\phi(t)$ is the remainder in the polynomial division of $f(t)$ by $g(t)$ (and hence, $\deg \phi(t)<n$). This is coherent with Proposition \ref{prop:trick}, because $f(t) \equiv \phi(t) \bmod g(t)$ and hence
    \[ f(C_g) \Lambda_n(t) = \phi(C_g) \Lambda_n(t) \equiv \phi(t)\Lambda_n(t) \bmod g(t) \equiv f(t)\Lambda_n(t) \bmod g(t).   \]
\end{remark}

For some special choices of $g(t)$, Companion Rings are known, and have been studied, by other names.  For example, if $g(t)=t^n-1$, the Companion Ring of $g(t)$ is the ring of \emph{circulant matrices} \cite{Davis}; if $g(t)=t^n+1$, we obtain \emph{skew-circulant matrices} \cite{Davis}; if $g(t)=t^n-k \ (k\in R)$, we get the ring of \emph{pseudo-circulant matrices} \cite{Vaidyanathan}; and if $g(t)=t^n$, the Companion Ring corresponds to \emph{lower triangular Toeplitz matrices} \cite{BP94}.

When $R$ is a field, and for any choice of the monic polynomial $g(t)$, the theory of companion matrices has been deeply studied and its relation to quotient rings is well known: see for example \cite{fastQR,DDPbounds,SteveFiedler} for both the general theory and the role of companion matrices (when $R$ is a subfield of $\mathbb{C}$) in numerical analysis, and \cite{NN16,NP15} for the link to quotient rings, as well as the references cited in these papers. The more general case where $R$ is a commutative ring, as  in this article, was studied in \cite{NW21}, where particular attention was given to the case where $R$ is an EDD.

\subsection{Cyclically presented groups}

Any finitely generated abelian group $A$ is isomorphic to a group of the form $A_0\oplus \Z^\beta$ where $A_0$ is a finite abelian group and the \emph{Betti number} (or \emph{torsion-free rank}) $\beta \geq 0$. 
Given a group presentation $P=\pres{x_0,\dots ,x_{n-1}}{r_0,\dots ,r_{m-1}}$ ($n,m\geq 1$), defining a group $G$, the \emph{relation matrix} of $P$ is the $m\times n$ integer matrix $M$ whose $(i,j)$ entry ($1\leq i\leq m, 1\leq j\leq n$) is the exponent sum of generator $x_{j-1}$ in relator $r_{i-1}$. If the rank of $M$ is $r$ and the non-zero invariant factors of the Smith Form of $M$ are  $s_1, \dots, s_{r}$ then the abelianization of $G$ is
\[ G^\mathrm{ab} \cong \Z_{s_1}\oplus \dots \oplus \Z_{s_r}\oplus \Z^{n-r}.\]
(See, for example, \cite[pp.\,146–149, Theorem 3.6]{MKS}.) Thus $\beta (G^\mathrm{ab}) =n-r$ and if $G^\mathrm{ab}=A_0\oplus \Z^\beta$ then the group order $|A_0|=\gamma_r(M)$. When it is clear from the context what presentation of a group is being considered, we may abuse terminology and refer to the relation matrix of a group.

The \emph{deficiency} of a group presentation is equal to the number of its generators minus the number of its relators. The deficiency of a presentation defining a group $G$ is bounded above by $\beta(G^\mathrm{ab})$, and so we define the deficiency of a group $G$ to be the maximum of the deficiencies of all presentations defining $G$. Deficiency zero groups play an important role both in group theory and in low dimensional topology. Groups of positive deficiency are infinite, whereas deficiency zero groups may be either finite or infinite. The fundamental group of every closed, connected, bounded 3-dimensional manifold has a deficiency zero presentation. Cyclic presentations and cyclically presented groups provide an important subclass of presentations and groups of deficiency zero that admit a cyclic symmetry. Specifically, a \emph{cyclic presentation} is a group presentation of the form
\[ P_n(w) = \pres{x_0,\dots ,x_{n-1}}{w(x_i,x_{i+1},\dots ,x_{i+n-1})\ (0\leq i<n)}\]
where $w=w(x_0, x_1, \dots , x_{n-1})$ is some fixed element of the free group with basis $\{x_0,\dots, x_{n-1}\}$ and the subscripts are taken $\bmod~n$, and the group $G_n(w)$ it defines is called a \emph{cyclically presented group}. If, for each $0\leq i<n$, the exponent sum of $x_i$ in $w(x_0, \dots , x_{n-1})$ is $a_i$ then, setting $R=\Z$, the relation matrix of $P_n(w)$ is the integer circulant matrix $M=f(C_g)^T$ where $g(t)=t^n-1$ and $f(t)=\sum_{i=0}^{n-1} a_it^i$ (and so the Smith forms of $M$ and $f(C_g)$ are equal).

\begin{example}
Consider the group $G$ defined by the cyclic presentation
\[   P_4(x_0x_1^{3}x_2^{-2})=\pres{x_0,x_1,x_2,x_3}{x_0x_1^{3}x_2^{-2},x_1x_2^{3}x_3^{-2},x_2x_3^{3}x_0^{-2},x_3x_0^{3}x_1^{-2}}.\]
Letting $f(t)=1+3t-2t^2$ and $g(t)=t^4-1$, then the relation matrix of this presentation is $M=f(C_g)^T$, where 
\[  f(C_g)=\begin{bmatrix}
    1 &0&-2&3\\
    3&1&0&-2\\
    -2&3&1&0\\
    0&-2&3&1\\
\end{bmatrix}, \quad \mathrm{and} \quad C_g =\begin{bmatrix}
    0&0&0&1\\
    1&0&0&0\\
    0&1&0&0\\
    0&0&1&0
\end{bmatrix}.    \]
The invariant factors of $f(C_g)$ are $1,1,3,48$, of which $3,48$ are the non-units, and hence $G^\mathrm{ab}\cong \Z_3\oplus \Z_{48}$.
\end{example}

Thus, results concerning the Smith forms of such matrices $f(C_g)$ provide information about the cyclically presented group $G_n(w)$. Not least, this can be used for testing non-isomorphism of groups: if $G_n(w)^\mathrm{ab}$, $G_{n'}(w')^\mathrm{ab}$ are not isomorphic then $G_n(w)$, $G_{n'}(w')$ are not isomorphic.

For a 3-manifold $M$, the first integral homology $H_1(M)$ is isomorphic to the abelianization of its fundamental group (see, for example, \cite[Theorem 2A.1]{Hatcher}). Thus, given a 3-manifold whose fundamental group has a cyclic presentation $P_n(w)$ with circulant relation matrix $f(C_g)^T$ ($g(t)=t^n-1$), the Smith form of $f(C_g)$ provides the homology of $M$. Many such families of groups and manifolds are described in \cite{CRS03}.

Moreover, partial information, in the form of lower bounds for the number of non-unit invariant factors of the relation matrix, can yield structural results concerning finiteness and solvability of the groups and  lower bounds for the Heegaard genus of the manifolds.  Following \cite{Wilson98}, we write $d(G)$ to denote the rank, or minimum number of generators, of a group $G$. The number of non-unit invariant factors of a relation matrix of a presentation of a group $G$ is equal to $d(G^\mathrm{ab})$, and so  (since $G^\mathrm{ab}$ is a quotient of $G$) provides a lower bound for $d(G)$. 

If $G$ is a finite group of deficiency zero then  $d(G^\mathrm{ab})\leq 3$ \cite[Theorem 9(ii)]{JohnsonWamsleyWright}.  A related result concerns solvable groups. These are groups that can be described in terms of abelian groups, through group extensions \cite[p.\,293]{MKS}.  If $G$ is a solvable group of deficiency zero then $d(G^\mathrm{ab})\leq 4$ \cite{Wilson91} (see \cite[Corollary 1.2]{WilsonNotes}). Thus lower bounds for $d(G^\mathrm{ab})$ can provide an effective tool for proving that groups of deficiency zero, such as cyclically presented groups, are infinite or non-solvable.

The \emph{Heegaard genus}, $g(M)$ of a  closed, connected, orientable 3-manifold is the minimum $g$ for which $M$ admits a Heegaard splitting of genus $g$ \cite{BoileauZieschang}; it is bounded below by the minimum number of generators $d(G)$ of the fundamental group $G=\pi_1(M)$, which in turn is bounded below by $d(G^\mathrm{ab})$. 

\section{Smith forms of matrices in Companion Rings}\label{sec:SNFmtxcomprings}

\subsection{Prior results}

\begin{theorem}[Non-zero invariant factors {\cite[Theorem A]{NW21}}]\label{thm:C}
Let $g(t) \in R[t]$ be monic of degree $n$, and let $f(t) \in R[t]$ where $R$ is an EDD. Suppose that $g(t) = G(t) z(t)$, $f(t) = F(t) z(t)$ where $z(t)$ is a monic common divisor of $f(t)$ and $g(t)$. Then  $f(C_g) \sim F(C_G) \oplus 0_{m\times m}$, where $m= \deg z(t)$. In particular, $F(C_G)$ has invariant factors $s_1,\dots,s_r$ if and only if $f(C_g)$ has invariant factors $s_1, \dots, s_r$ and $0, \dots, 0$ ($m$ times).
\end{theorem}

The following determinant formula is well known \cite[equation (1.1)]{NW21}:
\begin{equation}\label{eq:detfCg}\det f(C_g) = \prod_{\theta : g(\theta) = 0} f(\theta) =: \res(f,g).\end{equation}
 The immediate corollary of Theorem \ref{thm:C} below expresses the last non-zero determinantal divisor as the resultant of $F(t)$ and $G(t)$. This therefore generalizes the expression \eqref{eq:detfCg} to the case of singular matrices $f(C_g)$.

\begin{corollary}[Last non-zero determinantal divisor {\cite[Corollary B]{NW21}}]\label{cor:B}
In the notation of Theorem \ref{thm:C}, suppose that $z(t)$ is the monic \emph{greatest} common divisor of $f(t)=z(t)F(t)$ and $g(t)=z(t) G(t)$. Then the last non-zero determinantal divisor of $f(C_g)$ is
\[\gamma_r = \prod_{\theta : G(\theta) = 0} F(\theta) = \res(F,G).\]
\end{corollary}

On the other hand, the first determinantal divisor is given by the next result:

\begin{theorem}[First determinantal divisor {\cite[Lemma 7.1]{NW21}}]\label{thm:gamma1}
Let $f(t),g(t)\in R[t]$ where $R$ is an EDD and suppose $f(t)\equiv h(t) \bmod g(t)$ with $\mathrm{deg}(h(t))<\mathrm{deg} (g(t))$. Then $\gamma_1(f(C_g))=\cont(h)$. In particular, $\gamma_1(f(C_g))=1$ if and only if $h(t)$ is primitive.
\end{theorem}

The following results concern factorizations of $f(t)$ or $g(t)$:

\begin{theorem}[Factorizing $f(t)$ {\cite[Theorem 6.1]{NW21}}]
Let $f(t),g(t)\in R[t]$ where $R$ is an EDD. Let $f(t)=f_1(t) f_2(t)$ and suppose that $\res(f_1,g)$ and $\res(f_2,g)$ are coprime. Denote by $S,S_1,S_2$ the Smith forms of, respectively, $f(C_g)$, $f_1(C_g)$, $f_2(C_g)$. Then $S=S_1 S_2$.
\end{theorem}

\begin{corollary}[{\cite[Corollary 6.2]{NW21}}]\label{cor:factorf}
Let $f(t),g(t)\in R[t]$ where $R$ is an EDD. Let $f(t)=f_1(t) f_2(t)$ and suppose that $\res(f_2,g)$ is a unit of $R$. Then $f_1(C_g) \sim f(C_g)$.
\end{corollary}

\begin{theorem}[Factorizing $g(t)$ {\cite[Theorem 6.3]{NW21}}]\label{thm:factorgbetter} Let $f(t),g(t)\in R[t]$ where $R$ is an EDD. Let $g(t)=g_1(t) g_2(t)$ and suppose that $\res(f,g_1)$ and $\res(f,g_2)$ are coprime. Then $f(C_g) \sim f(C_{g_1}) \oplus f(C_{g_2})$.
\end{theorem}

\begin{corollary}[{\cite[Corollary 6.4]{NW21}}]\label{cor:coprime}
Let $f(t),g(t)\in R[t]$ where $R$ is an EDD. Let $g(t) = g_1(t) g_2(t)$ and suppose that $\res(f,g_2)$ is a unit of $R$. Then $f(C_g) \sim I_{\deg g_2(t)} \oplus f(C_{g_1}).$
\end{corollary}

\subsection{Second last determinantal divisor}

If $R$ is an integral domain and $f(t),g(t)\in R[t]$ then there exist $u(t),v(t)\in R[t]$ such that $f(t)u(t)+g(t)v(t)=\res (f,g)$ (see, for example, \cite[Lemma 7(1)]{Dresden}). If, in addition, $g(t)$ is monic and $f(t),g(t)$ are coprime polynomials, then $f(C_g)u(C_g)=\res (f,g) I_{\deg g(t)}$, so $u(C_g)= \res (f,g)  (f(C_g))^{-1} = \adj(f(C_g))$, and hence $\adj(f(C_g))=u(C_g)$ is an element of the Companion Ring of $g(t)$. Moreover, still assuming that $g(t)$ is monic, by Theorem \ref{thm:quotientremainder} there exists a unique $q(t)\in R[t]$ of degree less than $\deg g(t)$ such that $f(t)q(t)\equiv \res(f,g) \bmod g(t)$.

\begin{theorem}[Second last determinantal divisor]\label{thm:gamman-1}
Let $R$ be an EDD and let $f(t),g(t) \in R[t]$ be coprime integer polynomials, with $g(t)$ monic. Let $q(t)\in R[t]$ be the unique polynomial of degree less than $n = \deg g(t)$ such that $f(t) q(t) \equiv \res(f,g) \bmod g(t)$. Then, $\gamma_{n-1}(f(C_g)) = \cont (q(t))$. In particular, $\gamma_{n-1}(f(C_g))=1$ if and only if $q(t)$ is primitive.
\end{theorem}

\begin{proof}
By definition $\gamma_{n-1}$ is the GCD of the minors of $f(C_g)$ of order $n-1$. Such minors are the elements of $\adj(f(C_g))=q(C_g)$. Denote $q(t)=\sum_{k=0}^{n-1} q_kt^k$. Now, as noted in Section \ref{sec:companionrings}, since $C_g$ is a companion matrix, for each $0\leq k<n$, the bottom row of $C_g^{k}$ is equal to $e_{n-k}^T$ and so the bottom row of   $q(C_g)$ is $[q_{n-1}, \dots , q_1 , q_0]$. Thus $\gamma_{n-1}$ divides $\gcd (q_0,q_1,\dots,q_{n-1})=\cont(q)$. Moreover, the $(i,j)$-th entry of $q(C_g)$ is given by
\[ (q(C_g))_{ij}=\sum_{k=0}^{n-1}q_k(C_g^k)_{ij},\]
which is a linear combination, over $R$, of $q_0,\dots,q_{n-1}$. Thus each $(q(C_g))_{ij}$ is divisible by $\gcd(q_0,\dots,q_{n-1})=\cont(q)$, and so $\gamma_{n-1}$ is divisible by $\cont(q)$. Hence $\gamma_{n-1}=\cont(q)$, and the proof is complete.
\end{proof}

\begin{remark}
As observed in the proof of Theorem \ref{thm:gamman-1}, if $q(t)=\sum_{k=0}^{n-1} q_kt^k$ then the bottom row, $x^T$, say of $q(C_g)$ is equal to $[q_{n-1}, \dots , q_1 , q_0]$. That is, $x^T=e_n^Tq(C_g)=e_n^T\adj(f(C_g)) = e_n^T \res(f,g) (f(C_g))^{-1}$, or 
\begin{alignat}{1}
x^Tf(C_g) = e_n^T \res(f,g). \label{eq:linsystem}   
\end{alignat}
Thus computing $q(t)$ amounts to solving the linear system \eqref{eq:linsystem}.
\end{remark}

Combining Theorem \ref{thm:gamman-1} with Theorem \ref{thm:C} yields the following expression for the second last non-zero determinantal divisor of $f(C_g)$, in the more general setting where $f(t)$ and $g(t)$ are not coprime. 

\begin{corollary}\label{cor:gamman-1}
Let $f(t)=z(t) F(t)$ and $g(t)=z(t) G(t)$ be integer polynomials, where $g(t)$ is monic and $z(t)$ is the monic GCD of $f(t)$ and $g(t)$. Let $Q(t)$ be the unique polynomial of degree less than $r = \deg G(t)$ such that
\[ Q(t) F(t) \equiv \res(F,G) \bmod G(t).\]
Then $\gamma_{r-1}(f(C_g)) = \cont (Q(t))$. In particular, $\gamma_{r-1}(f(C_g))=1$ if and only if $Q(t)$ is primitive.
\end{corollary}

\subsection{Swap theorem}

Our main result in this section is Theorem \ref{thm:swaptheorem}
which, for monic polynomials $f(t),g(t)$, allows us to translate between $f(C_g)$ and $g(C_f)$. For this, we recall the concept of Horner shifts \cite{Horner0,Horner1}. Let $h(t)=\sum_{i=0}^m h_i t^i \in R[t]$ be a polynomial of degree $m$. The \emph{Horner shift of degree $0\leq k < m$} of the polynomial $h(t)$ is defined as
\begin{equation*}
    \sigma_k(h(t)) =\frac{h(t)-\sum_{i=0}^{m-k-1}h_i t^i}{t^{m-k}}.
\end{equation*}
For example, if $h(t)=t^3-2t+1 \in \Z[t]$, the associated Horner shifts are $\sigma_2(h(t))=t^2-2$, $\sigma_1(h(t))=t$, and $\sigma_0(h(t))=1$.

\begin{theorem}[Swap Theorem]\label{thm:swaptheorem}
Let $R$ be an EDD and let $f(t),g(t)\in R[t]$ be monic polynomials of degrees $m,n$, respectively, where $n\geq m$. Then $f(C_g)\sim I_{n-m}\oplus g(C_f)$.
\end{theorem}

We require the following technical result.

\begin{lemma}\label{lem:explicit}
Let $R$ be an EDD and let $g(t), f(t)\in R[t]$, where $g(t)$ is monic  of degree $n$. Let $\phi(t)$ be the unique polynomial of degree $m$, $0 \leq m<n$, and such that $\phi(t) \equiv f(t) \bmod g(t)$. Let $q(t),r(t)\in R[t]$ be the unique polynomials such that $t^{n-1}\phi(t)=q(t)g(t)+r(t)$, $\deg q(t)=m-1$ if $m>0$ (or $q(t)=0$ if $m=0$), and  $\deg r(t)<n$. Moreover, let
\[
\Psi(q(t))= \begin{bmatrix}
q(t) & \sigma_{m-2}(q(t)) & \cdots & \sigma_1(q(t)) &
\sigma_0(q(t)) & 0 & \cdots & 0\end{bmatrix}^T.
\]
(If $q(t)=0$, then $\Psi(q(t))=0$.)
Then
\[ f(C_g)\Lambda_n(t)= \phi(t) \Lambda_n(t)-g(t)\Psi(q(t)).\]
\end{lemma}

\begin{proof}
Note first that $f(C_g)=\phi(C_g)$. Moreover, by Theorem \ref{thm:quotientremainder}, the polynomials $q(t),r(t)\in R[t]$ defined in the statement are indeed unique, and it suffices to show
\begin{alignat}{1}
\phi(C_g)\Lambda_n(t)= \phi(t) \Lambda_n(t)-g(t)\Psi(q(t)).\label{eq:swaplemRTP}
\end{alignat} 
By Proposition \ref{prop:trick} $\phi(C_g)$ is the unique element of $R^{n\times n}$ such that  $\phi(C_g)\Lambda_n(t)\equiv \phi(t)\Lambda_n(t) \bmod g(t);$ hence, using also Theorem \ref{thm:quotientremainder}, the $(n-k)$th component of $\phi(C_g)\Lambda_n(t)$ is the remainder in the polynomial division of $\phi(t) t^k$ by $g(t)$, for all
$k=0,\dots,n-1$. This immediately establishes the top row of \eqref{eq:swaplemRTP} where $q(t)$ is the quotient in the same polynomial division (for $k=n-1$). Moreover, again by Theorem \ref{thm:quotientremainder}, for all $k=n-1,n-2,\dots,1$, there exist unique $a(t),b(t),c(t),d(t)\in R[t]$ with $\deg b(t),\deg d(t)<n$ such that 
\begin{alignat}{1}
    t^k\phi(t) &= g(t)a(t)+b(t),\label{eq:swaptk}\\
    t^{k-1}\phi(t) &= g(t)c(t)+d(t).\label{eq:swaptk-1}
\end{alignat}
We claim that $c(t)=(a(t)-a(0))/t$; by definition of Horner's shift, this proves \eqref{eq:swaplemRTP} by finite induction. To prove the claim, observe that \eqref{eq:swaptk},\eqref{eq:swaptk-1} imply
\[\frac{g(t)a(t)+b(t)}{t}=g(t)c(t)+d(t)\]
and so
\[g(t) \frac{a(t)-a(0)}{t}+\frac{a(0)g(t)+b(t)}{t}=g(t)c(t)+d(t).\]
\end{proof}

\begin{proof}[Proof of Theorem \ref{thm:swaptheorem}]
Suppose first $n>m$.

For any $1\leq k\leq m$ the matrix $C_g^k$ is of the form $\begin{bmatrix} *&*\\ I_{n-k}& 0\end{bmatrix}$, and hence we may partition $f(C_g) = \begin{bmatrix} A & B\\ X & C\end{bmatrix}$ where $B \in R^{m \times m}$ and $X$ is a unit upper triangular Toeplitz matrix. Now
\[ \begin{bmatrix} 0 & X^{-1}\\ I & -AX^{-1}\end{bmatrix}
    \begin{bmatrix} A & B\\ X & C\end{bmatrix}
    \begin{bmatrix} I & -XC^{-1}\\ 0&I \end{bmatrix}
=\begin{bmatrix} I&0\\0&B-AX^{-1}C \end{bmatrix}
\]
and so $f(C_g) \sim I_{n-m} \oplus ( B - A X^{-1} C  )$. Hence, it suffices to show that $g(C_f)\sim ( B - A X^{-1} C  )$.

By Lemma \ref{lem:explicit}, there exists a monic polynomial $q(t)=t^{m-1}+\sum_{i=0}^{m-2} q_it^i$ such that
\begin{alignat*}{1}
A t^m \Lambda_{n-m}(t) + B \Lambda_m(t) &= f(t) t^{n-m} \Lambda_{m}(t) - g(t) \Xi(t),\\
 X t^m \Lambda_{n-m}(t) + C \Lambda_m(t) &= f(t) \Lambda_{n-m}(t),
\end{alignat*}
where
\[
\Xi(t)=\begin{bmatrix}
    q(t) &
\sigma_{m-2}(q(t))&
\cdots &
\sigma_1(q(t)) &
\sigma_0(q(t))
\end{bmatrix}^T.
\]
Then 
\[ \left( A X^{-1} C  - B \right) \Lambda_m(t) \equiv  g(t) \Xi(t) \bmod f(t). \]
Introducing now the unit upper triangular Toeplitz matrix
\[U =\begin{bmatrix}
1 & q_{m-2} & \dots & q_1 & q_0\\
& 1 & q_{m-2} & \dots & q_1\\
& & \ddots & \ddots &\vdots\\
& & & 1&q_{m-2}\\
& & & & 1
\end{bmatrix},\]
we have $U \Lambda_m(t)=\Xi(t)$ and so
\[ U^{-1}\left( A X^{-1} C - B\right)\Lambda_m(t) \equiv g(t) \Lambda_m(t) \bmod f(t).\]
By Proposition \ref{prop:trick}, this implies $U^{-1}\left( A X^{-1} C - B\right)=g(C_f)$, and hence $g(C_f)\sim U g(C_f) =  A X^{-1} C - B $, which concludes the proof for the case $n>m$.

If $m=n$, then note that $A,C,X$ are empty matrices while $B=f(C_g)$. Moreover, in this case, $f(t) = g(t) + \phi(t)$ for some $\phi(t)$ such that $\deg \phi(t) := M < n$. By Lemma \ref{lem:explicit}, we then get
\[ -\phi(C_g) \Lambda_n(t) \equiv g(t) [\Lambda_n(t) + \Psi(q(t))] \bmod f(t)\]
(where $\Psi(q(t))$ is as defined in the lemma). Now let $V \in R^{n \times n}$ be defined by $V \Lambda_n(t) = \Lambda_n(t) + \Psi(q(t))$, and observe that $V$ must be unit upper triangular\footnote{In addition, one can prove that $V$ is also Toeplitz by the properties of $\Psi_q(t))$.}
, and hence unimodular, because both $\Lambda_n(t)$ and $\Lambda_n(t) + \Psi(q(t))$ are degree-graded vectors with monic components (see \cite{NN16,NNT17} for more details). Using also Proposition \ref{prop:trick}, it follows that
\[ -V^{-1} \phi(C_g) \Lambda_n(t) \equiv g(t) \Lambda_n(t) \bmod f(t), \]
implying $f(C_g)=\phi(C_g)=-V g(C_f)$; this implies the statement for the case $m=n$.
\end{proof}

\subsection{Composition theorem}

For $f(t),g(t),h(t)\in R[t]$, with $g(t),h(t)$ monic, the next result expresses the composition $(f \circ h) (C_{g \circ h})$ in terms of $f(C_g)$ and $\deg h(t)$. For this we recall that the \emph{Kronecker product} of $M$ and $N$, denoted $M \otimes N$, is the block matrix whose $(i,j)$ block entry is $M_{ij} N$ \cite{HJ91}. Note, in particular, that $I_k \otimes N$ is equal to the direct product of $k$ copies of $N$. For any pair of square matrices $M,N$ (possibly of different sizes) the Kronecker products $M\otimes N, N\otimes M$ are permutation similar. This is stated assuming that $R$ is (a subring of) a field in \cite[Corollary 4.3.10]{HJ91}, but since permutation matrices only contain $0$ or $1$ as their elements and have determinant $\pm 1$, the same proof is in fact valid for any commutative ring $R$.
\begin{theorem}[Composition Theorem]\label{thm:composition}
Let $R$ be an EDD and let $f(t),g(t),h(t)\in R[t]$ where $g(t),h(t)$ are monic. Then,
\[(f \circ h) (C_{g \circ h}) \sim_S f(C_g) \otimes I_{\deg h(t)} \sim_S I_{\deg h(t)} \otimes f(C_g)=\underbrace{f(C_g)\oplus\cdots \oplus f(C_g)}_{\deg h(t)}.\]
\end{theorem}

For the proof of Theorem \ref{thm:composition} we need the following:

\begin{lemma}\label{lem:composition}
Let $R$ be an EDD and let $g(t),h(t)\in R[t]$ be monic. Then $h (C_{g \circ h})\sim_S C_g \otimes I_{\deg h(t)}$.
\end{lemma}

\begin{proof}
Let $n=\deg g(t), m=\deg h(t)$ and write $h(t)=t^m+\sum_{i=0}^{m-1}h_it^i$. 
For notational simplicity set $H=h(C_{g\circ h})$ and $C=C_g  \otimes I_m$ and let
\[ A=\begin{bmatrix}
1 & h_{m-1} & \dots & & h_1\\
& 1 & h_{m-1} & &\vdots\\
& & \ddots & \ddots &\\
& & & \ddots &h_{m-1}\\
& & & & 1
\end{bmatrix}, 
\quad 
B=\begin{bmatrix}
h_0 &  &  & &\\
h_1& h_0 &  & &\\
& & \ddots & &\\
\vdots& & \ddots & \ddots &\\
h_{m-1}& & & h_1& h_0
\end{bmatrix} \in R^{m \times m},\]
\[V_k=\begin{bmatrix}
A & B & & & \\
& A & B & & \\
& & \ddots & \ddots & \\
& & & A & B \\
& & & & A 
\end{bmatrix} \in R^{mk \times mk},\]
\[U_k=V_k \oplus I_{nm-mk}\ (1\leq k\leq n), U=U_1 U_2 \cdots U_n.
\]
Since each $V_k$ is unit upper triangular we have $\det (U)=1$ and thus it suffices to show that $UH=CU$.
For each $1\leq k\leq n$ let
\begin{alignat*}{1}
\Theta_k(t)=\begin{bmatrix}
(A^{-1}(tI-B))^{k-1}A^{-1}\\
\vdots\\
A^{-1}(tI-B)A^{-1}\\
A^{-1}
\end{bmatrix}.
\end{alignat*}
Then $V_k\Theta_k = \begin{bmatrix} t\Theta_{k-1}\\ I_m \end{bmatrix}$, and so $U_n\Theta_n=V_n\Theta_n=\begin{bmatrix} t\Theta_{n-1}\\ I_m \end{bmatrix}$. Premultiplying by $U_{n-1},\dots ,U_1$ in turn gives $U \Theta_n= \begin{bmatrix} t^{n-1}I_m& t^{n-2}I_m& \cdots& tI_m &I_m\end{bmatrix}^T$,
or equivalently
\begin{alignat}{1}
\Theta_n(t)=U^{-1}(\Lambda_n(t)\otimes I_m).\label{eq:thetat}
\end{alignat}
Letting $S=RI_m$, the companion matrix of $g(t)I_m \in S[t]$ is $C_g\otimes I_m = C$. Hence, $C$ satisfies the elementwise (over $S$) congruence 
\[C (\Lambda_n(t) \otimes I_m) \equiv t (\Lambda_n(t) \otimes I_m) \bmod g(t) I_m,
\]
(see \cite{NN16,NNT17,NPer16} and the references therein). This implies, in particular, the elementwise (over $R$) congruence
\[C (\Lambda_n(t) \otimes I_m) \equiv t (\Lambda_n(t) \otimes I_m) \bmod g(t),
\]
and hence, by \eqref{eq:thetat}, we have
\begin{alignat*}{1}
CU\Theta_n(t)\equiv tU\Theta_n(t) \bmod g(t).
\end{alignat*}
Pre-multiplying by $U^{-1}$ and setting $X:=U^{-1}CU$ gives $X\Theta_n(t) \equiv t \Theta_n(t) \bmod g(t)$; then replacing $t$ by $h(t)$ gives 
\begin{equation}
X \Theta_n(h(t)) \equiv h(t) \Theta_n(h(t)) \bmod (g\circ h)(t).\label{eq:ineed}
\end{equation}
Now let
\[\sigma (h(t)) = \begin{bmatrix} \sigma_{m-1}(h(t)) & \cdots & \sigma_0(h(t))\end{bmatrix}^T. \]
By the definition of Horner shifts
\begin{equation}
 A\Lambda_m(t) = \sigma (h(t))\label{eq:ALambda=sigma}
\end{equation}
and it follows directly from the definition of $B$ that
\begin{equation}
(h(t)I-B) \Lambda_m(t) = t^{m}\sigma(h(t)).\label{eq:matrixLambda=tmsigma} 
\end{equation}
Then
\begin{alignat*}{1}
    \Theta_n(h(t))\sigma(h(t)) 
&=
\begin{bmatrix}
(A^{-1}(h(t)I-B))^{n-1}\Lambda_m(t)\\
\vdots\\
A^{-1}(h(t)I-B)\Lambda_m(t)\\
\Lambda_m(t)
\end{bmatrix}\quad \text{by}~\eqref{eq:ALambda=sigma}\\
&=
\begin{bmatrix}
t^{(n-1)m}\Lambda_m(t)\\
\vdots\\
t^{m}\Lambda_m(t)\\
\Lambda_m(t)
\end{bmatrix}\quad \text{using}~\eqref{eq:ALambda=sigma},\eqref{eq:matrixLambda=tmsigma}\\
&= \Lambda_{nm}(t).
\end{alignat*}
Post-multiplying \eqref{eq:ineed} by $\sigma(h(t))$ therefore gives
\[ X\Lambda_{nm}(t) \equiv h(t)\Lambda_{nm}(t) \bmod (g\circ h)(t)\]
and so Proposition \ref{prop:trick} implies $X=H$, as required.
\end{proof}

\begin{proof}[Proof of Theorem \ref{thm:composition}]
It suffices to prove the first similarity in the statement. Let $f(t)=\sum_i f_it^i$. By Lemma \ref{lem:composition} there exists a unimodular matrix $V$ such that $h(C_{g \circ h}) = V^{-1} (C_g \otimes I_{\deg h(t)}) V$. It follows that
\[
(f\circ h)(C_{g\circ h})
=\sum_i f_i V^{-1} (C_g \otimes I_{\deg h(t)})^i V 
= V^{-1} (f(C_g)\otimes I_{\deg h(t)}) V
\]
so $(f\circ h)(C_{g\circ h})\sim_S f(C_g)\otimes I_{\deg h(t)}$, as required.
\end{proof}

\subsection{The number of non-unit invariant factors}

Now let $R$ be a PID and for a fixed prime $p\in R$ we consider the quotient ring $R/\langle p \rangle$ of $R$ by the prime ideal generated by $p$. Note that $R/\langle p \rangle$ is a field and so the Smith form of a matrix over $R/\langle p \rangle$ has elements in $\{0,1\}$. (In particular, if $R=\Z$, then $R/\langle p \rangle$ is the finite field with $p$ elements, $\mathbb{F}_p$.) For $f(t)\in R[t]$ let $f_p(t):=[f(t) \bmod \langle p \rangle]$ denote the polynomial in $(R/\langle p \rangle)[t]$ such that $f(t)\equiv f_p(t)\bmod \langle p \rangle$. In this setting, the following theorem gives an expression for the number of non-unit invariant factors of $f(C_g)$, when $g(t)$ is monic.

\begin{theorem}[Number of non-unit invariant factors]\label{thm:genJO}
Let $R$ be a PID and let  $f(t),g(t) \in R[t]$ with $g(t)$ monic. Then $f(C_g)$ has precisely $\max_{p|\gamma_r} \deg(\gcd (f_p,g_p))$ non-unit invariant factors, where $\gamma_r$ is the last non-zero determinantal divisor of $f(C_g)$ and the maximum is taken over all primes $p \in R$ dividing $\gamma_r$.
\end{theorem}

\begin{proof}
Let $p \in R$ be a prime and, for $A \in R^{n \times n}$, define $\bar{A} \in (R/\langle p \rangle)^{n \times n}$  such that $A \equiv \bar{A} \bmod R/\langle p \rangle$. Fix an arbitrary minor $\mu$ in $A$ and the corresponding minor $\bar{\mu}$ in $\bar{A}$; then $\mu \equiv \bar{\mu} \bmod \langle p \rangle$. Let the $k$-th determinantal divisors of $A,\bar{A}$ be $\gamma_k,\bar{\gamma}_k$, respectively. Then $\bar{\gamma}_k=1$ if $p \nmid \gamma_k$ and $\bar{\gamma}_k=0$ if $p \mid \gamma_k$. Clearly $C_g \equiv C_{g_p} \bmod \langle p \rangle$, and hence $f(C_g) \equiv f_p(C_g) \equiv f_p(C_{g_p}) \bmod \langle p \rangle$. Using Theorem \ref{thm:C}, we conclude that $\deg(\gcd (f_p,g_p))=\ell_p$ if and only if $f_p (C_{g_p}) \sim I_{n-\ell_p} \oplus 0$ if and only if $p \mid \gamma_{n-\ell_p+1}(f(C_g))$ and $p \nmid \gamma_{n-\ell_p}(f(C_g))$. The statement follows because $p$ is arbitrary.

Finally, we note that it suffices to restrict to the set of primes that divide $\gamma_r$, for if $p \nmid \gamma_r$ then $\bar{\gamma}_k = 1$ for all $k \leq r$.
\end{proof}

\begin{remark}
    Even more generally, Theorem \ref{thm:genJO} can be stated for certain EDDs that are not PIDs. One example is when $R=\mathcal{A}(\Omega)$ is the ring of functions that are analytic over a connected open set $\Omega$. The crucial property that is needed in the proof is that, for all $x \in R$, $x$ is not a unit if and only if there is a prime that divides $x$. Note that this is not true of every EDD; for example, if $R=\mathbb{A}$ is the ring of algebraic integers, then $2 \in \mathbb{A}$ is not a unit but there is no prime that divides $2$. In fact, $\mathbb{A}$ contains no prime at all; see \cite[Remark 17]{N24}.
\end{remark}

We conclude this subsection by stating two corollaries of Theorem \ref{thm:genJO}, the second of which recovers a result of Johnson and Odoni \cite{JohnsonOdoni}:

\begin{corollary}\label{cor:genJO2}
    Let $R$ be a PID and let  $f(t),g(t) \in R[t]$ with $g(t)$ monic. Then, for all primes $p \in R$, $f(C_g)$ has at least $\deg(\gcd (f_p,g_p))$ non-unit invariant factors.
\end{corollary}

\begin{corollary}[{\cite[Proposition 4.1(ii)]{JohnsonOdoni}}]\label{cor:JohnsonOdoniFiniteCyclic} 
Let $f(t)\in \Z[t]$, $g(t)=t^n-1$ and suppose $\rho=\res (f,g)\neq 0$. Then $\gamma_{n-1}(f(C_g))=1$ if and only if $\gcd (f_p,g_p)$ is linear for every prime $p$ dividing $|\rho|$.
\end{corollary}

\begin{remark}
In Corollary \ref{cor:JohnsonOdoniFiniteCyclic}, the ring $\Z$ can be replaced by any PID. Combining this with Theorem \ref{thm:gamman-1} we obtain the following observation, which may be of independent interest: \em when $R$ is a PID, the polynomial $q$ of Theorem \ref{thm:gamman-1} is primitive if and only $\gcd (f_p,g_p)$ in $R/\langle p \rangle$ is linear for every prime dividing $|\rho|$. \em     
\end{remark}

\section{Applications}

In this section, we apply the results of Section \ref{sec:SNFmtxcomprings} to various problems in group theory and low-dimensional topology.

\subsection{The cocktail party graphs}\label{sec:circulantgraphs}

In the language of matrices in Companion Rings, \cite[Theorem 6.3]{WilliamsCirculantGraphs} can be stated as  Theorem \ref{thm:LAAThm6.3}, below. As described in \cite{WilliamsCirculantGraphs} this provides the Smith form of the adjacency matrix of the cocktail party graph  on $2n$ vertices (or hyperoctahedral graph or $(2n,n)$-Tur\'{a}n graph). 
 Theorem \ref{thm:LAAThm6.3} was proved in \cite{WilliamsCirculantGraphs} using Tietze transformations on cyclic presentations of groups. We now reprove it using the techniques developed in Section \ref{sec:SNFmtxcomprings}.

\begin{theorem}[{\cite[Theorem 6.3]{WilliamsCirculantGraphs}}]\label{thm:LAAThm6.3}
Let $m\geq 1$, $f(t)= (t^m+1)\sum_{i=0}^{m-2} t^i$, $g(t)=t^{2m}-1$. Then the invariant factors of $f(C_g)$ are $1$ ($m-1$ times), $m-1$ (1 time) and $0$ ($m$ times).
\end{theorem}

\begin{proof}
Let $z(t)=\gcd (f(t),g(t))$, $F(t)=f(t)/z(t)$, and $G(t)=g(t)/z(t)$. Then $z(t)=t^m+1$, $F(t)=H(t)-t^{m-1}$, $G(t)=(t-1)H(t)$ where $H(t)=\sum_{i=0}^{m-1} t^i$. Let $Q(t)=H(t)-(m-1)t$. Working $\bmod~G(t)$ we have $tH(t)\equiv H(t)$, so $H(t)^2=\sum_{i=0}^{m-1} t^iH(t)\equiv mH(t)$, and hence
\begin{alignat*}{1}
F(t)Q(t) 
&\equiv H(t)^2-t^{m-1}H(t)-(m-1)tH(t)+(m-1)t^m\\
&\equiv mH(t)-H(t)-(m-1)H(t)+(m-1)\\
&\equiv m-1=\res(F,G).
\end{alignat*}
Moreover, since $Q(t)$ is primitive, Corollary \ref{cor:gamman-1} implies that the second last determinantal divisor of $F(C_G)$ is trivial. The result then follows from Corollary \ref{cor:B}.
\end{proof}

\subsection{Fractional Fibonacci groups}\label{sec:FractionalFibonacci}

The \emph{fractional Fibonacci groups} $\mathcal{F}^{(k)}(n)$ ($k,n \geq 1$), introduced in \cite{Maclachlan}, are the cyclically presented groups $G_n(x_0x_1^kx_2^{-1})$, and they generalize Conway's Fibonacci groups $\mathcal{F}(n)=\mathcal{F}^{(1)}(n)$ \cite{Conway65}. For even $n=2m$ they are fundamental groups of closed, connected, orientable 3-dimensional \emph{fractional Fibonacci manifolds} $\bar{M}(k,m)$ \cite{Maclachlan,MaclachlanReid} and so $\mathcal{F}^{(k)}(2m)^\mathrm{ab}$ provides the first integral homology of $\bar{M}(k,m)$. The contrasting case, $n$ odd, is investigated in \cite{ChinyereWilliamsFF}.

The relation matrix of $\mathcal{F}^{(k)}(n)$ is the circulant matrix $f(C_g)^T$ where $g(t)=t^n-1$ and $f(t)=t^2-kt-1$, which has roots $\lambda_\pm=(k\pm \sqrt{k^2+4})/2$. As in \cite{Maclachlan}, we define the \emph{fractional Fibonacci numbers}
\begin{alignat}{1}
F_0^k&=0, F_1^k=1, F_{j+2}^k=kF_{j+1}^k+F_{j}^k\ (j\geq 0) \label{eq:FFseq}
\end{alignat}
(which are the classical Fibonacci numbers $F_j$ in the case $k=1$) and it follows that
\begin{alignat}{1}
    F_n^k=\frac{\lambda_+^n-\lambda_-^n}{\lambda_+-\lambda_-}.\label{eq:FFnumber}
\end{alignat} 
Maclachlan \cite[Section 3]{Maclachlan} observed that
\begin{alignat}{1}
|\mathcal{F}^{(k)}(n)^\mathrm{ab}|=F_{n+1}^k+F_{n-1}^k-1-(-1)^n. \label{eq:FFabOrder}
\end{alignat}
In Theorem \ref{thm:FF} we give a formula for the structure of  $\mathcal{F}^{(k)}(n)^\mathrm{ab}$ that involves GCDs of expressions in the numbers $F_j^k$. By simplifying these GCDs we obtain an alternative formula for $\mathcal{F}^{(k)}(n)^\mathrm{ab}$ in Corollary \ref{cor:FFab}. In the case $k=1$ this coincides with the formula for $\mathcal{F}(n)^\mathrm{ab}$ given in \cite{Kim78}, and in the case $n$ even it coincides with \cite[Lemma 1]{MaclachlanReid}, which was stated without proof. One may infer insights about $\mathcal{F}^{(k)}(n)$ from Corollary \ref{cor:FFab} that are not evident from Theorem \ref{thm:FF}. For example, an expectation implicit in \cite[Theorem 6.2(a)]{ChinyereWilliamsFF}, that if $k$ is odd and $n\equiv 3\bmod 6$ then $\mathcal{F}^{(k)}(n)^\mathrm{ab}\cong (Q_8\times \Z_{(F_{n+1}^k+F_{n-1}^k)/4})^\mathrm{ab}$, where $Q_8$ denotes the quaternion group.

\begin{theorem}\label{thm:FF}
   Let  $n,k\geq 1$. Then  $\mathcal{F}^{(k)}(n)^\mathrm{ab} \cong \Z_\alpha \oplus \Z_\beta$, where $$\alpha=\gcd(F_n^k,F_{n-1}^k-1) \quad and \quad   \beta = 
   \frac{F_{n+1}^k+F_{n-1}^k-1-(-1)^n}{\alpha}.$$
\end{theorem}

\begin{proof}
 As noted above, the relation matrix of $\mathcal{F}^{(k)}(n)$ is $f(C_g)^T$ where $f(t)=t^2-kt-1$, $g(t)=t^n-1$. By Theorem \ref{thm:swaptheorem} we have $f(C_g)\sim I\oplus g(C_f)$, so it suffices to consider $g(C_f)$. Define $h(t)=F_n^k t+(F_{n-1}^k-1)$. We claim $g(t)\equiv h(t) \bmod f(t)$ for all $n\geq 1$. The case $n=1$ is immediate. With the inductive hypothesis $t^n-1 \equiv F_n^k t+(F_{n-1}^k-1) \bmod f(t)$, and working $\bmod f(t)$, we have
\begin{alignat*}{1}
    t^{n+1}-1 
    &\equiv t(t^n-1) + (t-1)\\
    &\equiv F_n^k t^2 + (F_{n-1}^k-1)t + (t-1)\\
    &\equiv F_n^k (kt+1) + F_{n-1}^kt -1\\
    &\equiv (k F_n^k + F_{n-1}^k)t + F_n^k-1\\
    &\equiv F_{n+1}^kt + (F_n^k-1),
\end{alignat*}
proving the claim.

By Theorem \ref{thm:gamma1} $\gamma_1(g(C_f))=\cont(h)=\gcd(F_n^k,F_{n-1}^k-1)$ and, by \eqref{eq:FFabOrder}, $\gamma_2(g(C_f))=\res (f,g)=F_{n+1}^k+F_{n-1}^k-1-(-1)^n$. The result follows.
\end{proof}

\begin{corollary}\label{cor:FFab}
Let $n,k\geq 1$. Then
\[\mathcal{F}^{(k)}(n)^\mathrm{ab}\cong
\begin{cases}
\Z_{F_{n+1}^k+F_{n-1}^k} & \text{if}~n\equiv 1,5,7,11\bmod 12,\\
\Z_{\gcd (k+1,2)}\oplus \Z_{\frac{F_{n+1}^k+F_{n-1}^k}{\gcd (k+1,2)}} & \text{if}~n\equiv 3,9\bmod 12,\\
\Z_{F_{n/2+1}^k+F_{n/2-1}^k}\oplus  \Z_{F_{n/2+1}^k+F_{n/2-1}^k} & \text{if}~n\equiv 2,6,10 \bmod 12,\\
\Z_{\gcd (k,2)F_{n/2}^k}\oplus \Z_{\frac{(k^2+4)F_{n/2}^k}{\gcd(k,2)}} & \text{if}~n\equiv 0,4,8\bmod 12.
\end{cases}
\]
\end{corollary}

\begin{proof}[Sketch proof of Corollary \ref{cor:FFab}]
Let $\alpha,\beta$ be as given in Theorem \ref{thm:FF}. An inductive argument,  involving recurrence relations, the definition \eqref{eq:FFseq} and the formula \eqref{eq:FFnumber} shows that for all odd $j$, $1\leq j\leq n/2$
\begin{alignat}{1}
\alpha=\gcd( F_{n-j}^k-F_j^k, F_{n-(j+1)}^k+F_{j+1}^k).\label{eq:fn-j}
\end{alignat}
Suppose first that $n$ is odd. If $(n-1)/2$ is odd (resp.\,$(n-1)/2$ is even) then substituting $j=(n-1)/2$ (resp.\,$(n+1)/2$) into \eqref{eq:fn-j} and simplifying gives $\alpha=\gcd (F_{(n+1)/2}^k-F_{(n-1)/2}^k, 2F_{(n+1)/2}^k)$. Then $F_{(n+1)/2}^k-F_{(n-1)/2}^k$ is even if and only if $k$ is odd and $n\equiv 3,9 \bmod 12$, in which case $\alpha=2$; otherwise $\alpha=1$. Therefore, if $n\equiv 1,5,7,11 \bmod 12$ then $g=1$ and if $n\equiv 3,9\bmod 12$ then $\alpha=\gcd(k,2)$, and the result follows.

Consider now the case $n$ even. If $n\equiv 2\bmod 4$ then substituting $j=n/2$ into \eqref{eq:fn-j} gives $\alpha=F_{n/2-1}^k+F_{n/2+1}^k$, and if $n\equiv 0\bmod 4$ then substituting $j=n/2-1$ into \eqref{eq:fn-j} and using \eqref{eq:FFnumber} gives $\alpha=\gcd (k,2)F_{n/2}^k$. The following identity can be confirmed by, on each side, expressing each term $F_j^k$ according the formula \eqref{eq:FFnumber} and simplifying using  $\lambda_+\lambda_-=-1$:
    \[ F_{n+1}^k+F_{n-1}^k-2=\begin{cases}
        (F_{n/2-1}^k+F_{n/2+1}^k)^2 & \text{if}~n\equiv 2\bmod 4,\\
        (k^2+4)(F_{n/2}^k)^2 & \text{if}~n \equiv 0\bmod 4,
    \end{cases}\]
and the value of $\beta$ follows.
\end{proof}

\subsection{Periodic generalized Neuwirth groups }\label{sec:neuwirth}

Let $n\geq 1$ and let $\alpha,\beta$ be coprime integers with $\alpha\geq 2\beta$, and let $l\geq 1$. The \emph{periodic generalized Neuwirth groups} are the groups
\[\Gamma_n((\alpha,\beta);\ell) = G_n( (x_0^\beta x_1^\beta\dots x_{n-1}^\beta)^\ell x_{n-1}^{-\alpha})
\]
and they are fundamental groups of closed, connected, orientable, 3-manifolds  \emph{periodic generalized Neuwirth manifolds} $M_n((\alpha,\beta);\ell)$ \cite[Section 3]{SpaggiariNeuwirth}. They form a subclass of the \emph{generalized Neuwirth groups} $\Gamma ((\alpha_1,\beta_1),\dots , (\alpha_n,\beta_n);\ell)$ defined in \cite{SpaggiariNeuwirth}. The groups $\Gamma_n((k+1,1);1)$ are the cyclically presented groups $G_n(x_0x_1\dots x_{n-2}x_{n-1}^{-k})$ considered in \cite{SVNeuwirth} and the groups $\mathcal{F}(n-1,n-1,k,n)$ of \cite{JohnsonOdoni}; in the case $k=1$ they coincide with the \em Neuwirth groups \em considered in \cite{Neuwirth68} and the generalized Fibonacci groups $\mathcal{H}(n-1,n,1)$ of \cite{CampbellRobertson75}. The first integral homology of $M_n((\alpha,\beta);\ell)$ is stated without proof in \cite[Lemma 3.2]{SpaggiariNeuwirth}. In Corollary \ref{cor:periodicgenNeuwirth} we prove this fact as an immediate corollary of the main result of this section, Theorem \ref{thm:circaaabbb}.

 The relation matrix for $\Gamma_n((\alpha,\beta);\ell)$ is $f(C_g)^T$ where
$g(t)=t^n-1$ and
\[f(t)=b\sum_{i=0}^{s-1}t^i+a\sum_{i=s}^{n-1}t^i=a\frac{t^n-1}{t-1}+(b-a)\frac{t^s-1}{t-1}.\]
For the above choices of $f(t)$ and $g(t)$, these matrices $f(C_g)$ have received some attention in the literature. For example, \cite[Exercise 27, p.81]{Davis} asks to obtain their determinant. In the case $a=0,b=1$, the Smith form of $f(C_g)$ was obtained in \cite[p. 184]{Rushanan} and \cite[Section 3]{WilliamsCirculantGraphs}, and was shown to be non-singular if and only if $\gcd(n,s)=1$, in \cite[Theorem 1]{BustomiBarra}. In Theorem \ref{thm:circaaabbb} we calculate the Smith form of $f(C_g)$ when $\gcd(n,s)=1$. Lemma \ref{lem:circ111000inverse}, which may be of independent interest, calculates the inverse of $f(C_g)$ when $a=0,b=1$. Hence, the inverse in the general case may be readily computed (as in the proof of Theorem \ref{thm:circaaabbb}).

\begin{lemma}\label{lem:circ111000inverse}
Let $n > s > 1$ where $\gcd(n,s)=1$ and let $L=f(C_g)$ where $f(t)=\sum_{i=0}^{s-1}t^i$ and $g(t)=t^n-1$. Moreover, let $v$ satisfy $0\leq v<s$ and $vn\equiv 1\bmod s$, and let $0\leq r<s$ satisfy $r\equiv n \bmod s$. Then $L^{-1}=q(C_g)$, with $q(t)=\sum_{i=0}^{n-1} q_i t^i$ and where, for $0\leq i<s$, the values of $q_i$ are as follows:
\[ q_{rj \bmod s} = \begin{cases}
    v/s-1 & 1\leq j\leq v,\\
    v/s & v+1\leq j\leq s,
\end{cases}\]
and, for $s\leq i<n$, $q_i=q_{i \bmod s}$.
\end{lemma}
\begin{proof}
Since $L$ is a circulant matrix, so is its inverse. That is, $L^{-1}=q(C_g)$, where $q(t)=\sum_{i=0}^{n-1} q_it^i$ for some $q_0,\ldots,q_{n-1}$, which has first column $L^{-1}e_1=\begin{bmatrix} q_0,\dots , q_{n-1} \end{bmatrix}^T=:q$, say. Then $Lq=e_1$, which corresponds to $n$ equations in the $n$ variables $q_0,\dots ,q_{n-1}$. The first such equation is
\[  q_0+q_{n-1}+\dots + q_{n-(s-1)} =1. \]
In addition, for each $2\leq j\leq n$, subtracting the $j$-th equation from the $(j-1)$-th equation yields the $(n-1)$ equations
\begin{alignat}{1}
    -q_1 + q_{n-s+1}&=1,\nonumber \\
    -q_j + q_{n-s+j}&=0 \quad (2\leq j <s),    \label{eq:Lj1}\\
    -q_{j} + q_{j-s}&=0 \quad (s \leq j <n).    \label{eq:Lj2}
\end{alignat}
Substituting $j=n-s,\dots , n-1$ into \eqref{eq:Lj1} and \eqref{eq:Lj2} implies that for $0\leq i,j <n$, if $i\equiv j\bmod s$ then $q_i=q_j$. Thus we may eliminate $q_s,\dots , q_{n-1}$ to leave $s$ variables $q_0,\dots , q_{s-1}$ and $s$ equations:
\begin{alignat}{1}
    q_0-q_r+\sum_{i=0}^{s-1} q_i &=1,\label{eq:LL0}\\
    -q_1 + q_{r+1}&=1,\label{eq:LL1}\\
    q_j &= q_{r+j} \quad (2\leq j\leq s-1).\label{eq:LLj}
\end{alignat}
For each $0\leq j<s$ define now $p_j:=q_{rj \bmod s}$. Then, since $vr\equiv 1\bmod s$, for each $0\leq j<s$ we have $q_j=p_{vj \bmod s}$. In particular $q_0=p_0$, $q_1=p_v$, $q_r=p_1$, $q_{r+1}=p_{v+1}$. For each $0\leq j<s$, define $t= vj \bmod s$, $0 \leq t < s$.
Then equations \eqref{eq:LL0},\eqref{eq:LL1},\eqref{eq:LLj} become
\begin{alignat}{1}
    p_0-p_1+\sum_{i=0}^{s-1} p_i &=1,\label{eq:pLL0}\\
    -p_v + p_{v+1}&=1,\label{eq:pLL1}\\
    p_t &= p_{t+1} \quad (0\leq t\leq s-1, t\neq 0,v).\nonumber%
\end{alignat}
Therefore $p_1=p_2=\dots = p_v$, $p_{v+1}=p_{v+2}=\dots =p_{s-1}=p_0$ and so equations \eqref{eq:pLL0},\eqref{eq:pLL1} become
\begin{alignat*}{1}
    (v-1)p_1+(s-v+1)p_0 &=1,\\
    -p_1+p_0&=1.%
\end{alignat*}
Solving gives $p_0=v/s, p_1=v/s-1$. Therefore $q_{rj \bmod s}=p_j=v/s-1$ if $1\leq j \leq v$ and $q_{rj \bmod s}=p_j=v/s$ if $v+1\leq j\leq s$.
\end{proof}

\begin{theorem}\label{thm:circaaabbb}
Let $n>s\geq 1$, $g(t)=t^n-1$,
\[f(t)=b\sum_{i=0}^{s-1}t^i+a\sum_{i=s}^{n-1}t^i,\]
where $a,b \in \Z$, and $\gcd(n,s)=1$. Define $k:=\frac{|a(n-s)+sb|}{\gcd(a,b)}$. Then the invariant factors of $f(C_g)$ are
    \[  \gcd (a,b), \underbrace{|a-b|,\dots,|a-b|}_{n-2 \ \mathrm{times}}, k \cdot |a-b|. \]
\end{theorem}

\begin{proof}
If $a=b$ then $f(C_g)=aee^T$, and the statement is readily obtained, so assume $a\neq b$. If $\gcd(a,b)>1$ then, letting $h(t) = \frac{f(t)}{\gcd(a,b)}$, we have $f(C_g) = \gcd(a,b) h(C_g)$, and thus the invariant factors of $f(C_g)$ are equal to those of $h(C_g)$ times $\gcd(a,b)$. Thus we may assume $\gcd(a,b)=1$.

Therefore $k=f(1)=a(n-s)+bs$ and by Theorem \ref{thm:gamma1}, $\gamma_1(f(C_g))=1$. Now, $\gamma_2$ is the GCD of all minors of order $2$. The $2\times 2$ submatrices of $f(C_g)$ are of the form $\begin{bmatrix} x&y\\ z&t\end{bmatrix}$ where $x,y,z,t\in\{a,b\}$ and so have determinant $0,\pm (a^2-b^2), \pm a(a-b), \pm b(a-b)$ (all of which do arise), and so $\gamma_2=|a-b| (0,a,b+a)=|a-b|$.

We have 
\begin{alignat*}{1}
\gamma_n
&=\res(f,g)=f(1)\cdot \res\left((b-a)\sum_{i=0}^{s-1}t^i+a\sum_{i=0}^{n-1} t^i,\sum_{i=0}^{n-1} t^i\right)\\
&\qquad =k\cdot \res\left((b-a)\sum_{i=0}^{s-1}t^i,\sum_{i=0}^{n-1} t^i\right)=k\cdot |a-b|^{n-1}.
\end{alignat*}

We now consider $\gamma_{n-1}$. Assume first $k=0$. Then $n-s$ divides $b$ and
\[ 1 = \gcd(a,b)=\gcd\left(-s \frac{b}{n-s},b\right)=\frac{|b|}{n-s} \gcd(s,n-s) \]
and so $|b|=n-s$, $|a|=s$ and $|a-b|=n$. Then $f(1)=k=0$ so $t-1$ divides $f(t),g(t)$ and Theorem \ref{thm:C} implies 
\[ \gamma_{n-1} = \res\left(\frac{f(t)}{t-1},\frac{g(t)}{t-1}\right) =  \frac{|a-b|^{n-1}}{n} \res \left( \frac{t^s-1}{t-1}, \frac{t^n-1}{t-1}   \right) = \frac{|a-b|^{n-1}}{n} = |a-b|^{n-2}.\] 

If instead $k\neq 0$, then $f(C_g)$ is invertible over $\mathbb{Q}$.  Let $L=\theta(C_g)$, where $\theta(t)=\sum_{i=0}^{s-1} t^i$. Then $f(C_g)=(b-a)L+aee^T$. Using the Sherman-Morrison formula \cite{SM50}, we have
\begin{alignat*}{1}
f(C_g)^{-1} 
&= \frac{L^{-1}}{b-a} - \frac{a L^{-1} ee^T L^{-1}}{(b-a)(b-a+a e^T L^{-1} e)}\\
&= \frac{L^{-1}}{b-a} - \frac{a ee^T}{ks(b-a)}.
\end{alignat*}
Therefore
\begin{alignat*}{1}
\adj(f(C_g)
&=\det (f(C_g))f(C_g)^{-1}\\
&=(b-a)^{n-2} \left(k L^{-1} - \frac{a}{s} ee^T \right).
\end{alignat*}
If $s=1$, then $L=L^{-1}=I$; otherwise, $L^{-1}$ is given by Lemma \ref{lem:circ111000inverse}. In either case, $(b-a)^{2-n}\adj(f(C_g)$ has precisely two distinct elements. If $s=1$ then these elements are $k-a$ and $-a$ (which are coprime).  If $s\neq 1$ then they are equal to $N=(kv-a)/s$ and $N-k$, where $v$ is as defined in Lemma \ref{lem:circ111000inverse}. Now 
\[ \gcd(Ns,k)=\gcd(kv-a,k)=\gcd(a,k)=\gcd(a,a(n-s),b)=\gcd(a,b)=1\]
so $\gcd(N,N-k)=\gcd(N,k)=1$.
Thus the GCD of the entries of $\adj(f(C_g)$ is $(b-a)^{n-2}$ and so again $\gamma_{n-1}=|a-b|^{n-2}$.

Now let $s_i$ denote the $i$-th invariant factor of $f(C_g)$. Then: $s_1=\gamma_1=1$; $s_n=\gamma_n/\gamma_{n-1}=k|a-b|$; $\gamma_2=|a-b|$ divides $s_2$, and hence for all $2\leq i \leq n$, we have $s_i=\mu_i|a-b|$ for some positive integers $\mu_i$. But $s_1\dots s_{n-1}=\gamma_{n-1}=|a-b|^{n-2}$ so $\mu_2=\dots =\mu_{n-1}=1$, or equivalently $s_2=\dots=s_{n-1}=|a-b|$, as required. 
\end{proof}

\begin{corollary}[{\cite[Lemma 3.2]{SpaggiariNeuwirth}}]\label{cor:periodicgenNeuwirth}
Let $n\geq 2$, $l,\alpha,\beta \geq 1$ where $(\alpha,\beta)=1$. Then the first integral homology of the manifold $M_n((\alpha,\beta);\ell)$ is isomorphic to $\Z_{\alpha}^{n-2}\oplus \Z_{\alpha |nl\beta -\alpha|}$.
\end{corollary}

\subsection{Generalized Fibonacci groups}\label{sec:Hrns}

The \emph{generalized Fibonacci groups} $\mathcal{H}(r,n,s)$ ($r,s,n\geq 1$), introduced in \cite{CampbellRobertson75}, are the cyclically presented groups $G_n(w)$ where
\[ w=x_0x_1\dots x_{r-1}(x_rx_{r+1}\dots x_{r+s-1})^{-1}.\]
They generalize the Fibonacci groups $\mathcal{F}(n)=\mathcal{H}(2,n,1) \cong \mathcal{H}(1,n,2)$. Without loss of generality we may assume $s\geq r$. The relation matrix of $\mathcal{H}(r,n,s)$ is the circulant matrix $f^{r,s}(C_g)^T$ where $g(t)=t^n-1$ and
\[f^{r,s}(t)  = 1+t+\dots +t^{r-1} -t^r(1+t+\dots +t^{s-1}).\]
Let $d =\gcd (r, n, s), R=r/d, N=n/d, S=s/d$. The groups $\mathcal{H}(r,n,s)$ with trivial abelianization were classified in \cite[Theorem A]{BainsonChinyere}. In \cite[Corollary 3.2]{NW22}, Corollary \ref{cor:B} was used to show that that if $s>r$ then $\mathcal{H}(r,n,s)^\mathrm{ab}\cong A_0\oplus \Z^{d-1}$ where  $|A_0|=|\mathcal{H}(R,N,S)^\mathrm{ab}|^d/|S-R|^{d-1}$ and that if $s=r$ then $\mathcal{H}(r, n, s)^\mathrm{ab}\cong A_0\oplus \Z^d$ where $|A_0|=N^{d-1}$.  In Theorem \ref{thm:Hrns}, for the case $S-R=1$, we express $\mathcal{H}(r,n,s)^\mathrm{ab}$ in terms of $\mathcal{H}(R,N,S)^\mathrm{ab}$, and in Theorem \ref{thm:H(r,n,r)} we calculate the structure of $\mathcal{H}(r,n,r)^\mathrm{ab}$.

\begin{theorem}\label{thm:Hrns}
Let $n\geq 2$, $s>r\geq 1$, $d=\gcd(r,n,s)$, $R=r/d, N=n/d, S=s/d$, and suppose $S-R=1$. Then 
\[\mathcal{H}(r,n,s)^\mathrm{ab}\cong \left( \mathcal{H}(R,N,S)^{\mathrm{ab}}\right)^d \oplus \Z^{d-1}.\]
\end{theorem}

\begin{proof}
As described above, the relation matrix of $\mathcal{H}(r,n,s)$ is $f(C_g)^T$, where $f(t)=f^{r,s}(t)$ and $g(t)=t^n-1$. Let $z(t)=\gcd (f(t),g(t))$, $F(t)=f(t)/z(t)$, and $G(t)=g(t)/z(t)$. Then if $h(t)=t^d$, $k^{(n)}(t)=\sum_{i=0}^{n-1} t^i$  we have $g(t)=(t-1)k^{(n)}(t)$ and (as shown in \cite[Proof of Theorem 3.1(a)]{NW22}, \cite[Proof of Theorem C]{WilliamsLOG}), $z(t) =\sum_{i=0}^{d-1} t^i$, $F(t)=(f^{R,S}\circ h)(t)$, $G(t)=(t-1) (k^{(N)}\circ h)(t)$.

By Theorem \ref{thm:C} the invariant factors of $f(C_g)$ are the invariant factors of $F(C_G)$ together with $0$ ($d-1$ times). Now $F(1)= R-S=-1$ so 
\[\gcd\left( \res(F(t),t-1), \res(F(t), k^{(N)}(t)\circ h(t)) \right)=1\]
and by Theorem \ref{thm:factorgbetter} $F(C_G) \sim F(C_{k^{(N)}\circ h}) \oplus F(C_{t-1})$. Then Theorem \ref{thm:composition} implies 
\begin{alignat*}{1}
    F(C_{k^{(N)}\circ h})&=(f^{R,S}\circ h) ( C_{k^{(N)}\circ h} )\\
    &\sim I_{\mathrm{deg}(h)} \otimes f^{R,S}(C_{k^{(N)}})\\
    &= I_{\mathrm{deg}(h)} \otimes f^{R,S}(C_{t^N-1})\quad \text{by Corollary \ref{cor:factorf}}\\
    &= \underbrace{f^{R,S}(C_{t^N-1})\oplus \cdots \oplus f^{R,S}(C_{t^N-1})}_{d \ \mathrm{times}}.
\end{alignat*}
Since $f^{R,S}(C_{t^N-1})^T$  is the relation matrix of $\mathcal{H}(R,N,S)$, the result follows.
\end{proof}

\begin{theorem}\label{thm:H(r,n,r)}
Let $n\geq 2, r\geq 1$ and let $d=\gcd(n,r)$, $N=n/d$. Then $\mathcal{H}(r,n,r)^\mathrm{ab} \cong \Z_N^{d-1}\oplus \Z^{d}$.
\end{theorem}

\begin{proof}
 The relation matrix of $\mathcal{H}(r,n,r)$ is $f(C_g)^T$ where $f(t)=f^{r,r}(t)$ and $g(t)=t^n-1$. Let $z(t)=\gcd (f(t),g(t))$, $F(t)=f(t)/z(t)$, and $G(t)=g(t)/z(t)$. Then by \cite[Proof of Theorem 3.1(b)]{NW22} we have $z(t)=t^d-1$, $G(t)=\sum_{i=0}^{N-1} t^{id}$, and $F(t)= F_0(t)^2F_1(t)$ where $F_0(t)=\sum_{i=0}^{R-1} t^{id}$, $F_1(t)=\sum_{i=0}^{d-1} t^{i}$, and, moreover, $\res(F_0,G)=1$. By Corollary \ref{cor:factorf} and Theorem \ref{thm:swaptheorem} $F(C_G) \sim F_1(C_G) \sim I_{n-2d+1} \oplus G(C_{F_1})$. Using the Cayley-Hamilton theorem, we have $C_{F_1}^d = I_{d-1}$ so $G(C_{F_1}) =N I_{d-1}$. Thus $f(C_g)\sim I_{n-2d+1} \oplus N I_{d-1} \oplus 0_{d\times d}$ and so $\mathcal{H}(r,n,r)^\text{ab}\cong \Z_N^{d-1}\oplus \Z^d$, as required.
\end{proof}

\subsection{Cyclically presented groups with length three positive relators}\label{sec:positivelength3}

 The cyclically presented groups with length three positive relators are the groups $G_n(x_0x_kx_l)$ ($0\leq k,l<n$), and have been studied in \cite{CRS05,EdjvetWilliams,MohamedWilliams,ChinyereEdjvetWilliams}. In \cite[Theorem 9.1]{MohamedWilliams}, the authors identified various pairs of groups $G_n(x_0x_{k_1}x_{l_1})$, $G_n(x_0x_{k_2}x_{l_2})$ that may have isomorphic abelianizations for all $n$, and they confirmed computationally that this is the case for $n\leq 1000$. In Theorem \ref{thm:MWCase(1)} we prove that, for the first such pair, the abelianizations are indeed isomorphic.
 
\begin{theorem}\label{thm:MWCase(1)}
If $16|n$ then $G_n(x_0x_1x_{n/2})^\mathrm{ab}\cong G_n(x_0x_1x_{n/4})^\mathrm{ab}\cong \Z_{2^{n/2}-1}$.
\end{theorem}

\begin{proof}
 The relation matrices of $G_n(x_0x_1x_{n/2}), G_n(x_0x_1x_{n/4})$ are $f_1(C_g)^T, f_2(C_g)^T$ where $g(t)=t^n-1$ and $f_1(t)=1+t+t^{n/2}, f_2(t)=1+t+t^{n/4}$, respectively. Thus it suffices to show $f_i(C_g)\sim I_{n-1}\oplus (2^{n/2}-1)$ for $i=1,2$. We shall write $g(t)=g_1(t)g_2(t)=h_1(t)h_2(t)g_2(t)$ where $g_1(t)=t^{n/2}-1$, $g_2(t)=t^{n/2}+1$, $h_1(t)=t^{n/4}-1$, $h_2(t)=t^{n/4}+1$.

Consider first $f_1$. The resultant $\res(f_1,g_2)=1$, so $f_1(C_{g_2})\sim I_{n/2}$, and hence by Corollary \ref{cor:coprime} we have $f_1(C_{g})\sim I_{n/2}\oplus f_1(C_{g_1})$. Moreover, $f_1(C_{g_1})=\phi(C_{g_1})$ where $\phi(t)=2+t$, so by Theorem \ref{thm:swaptheorem}
\[ f_1(C_{g_1}) \sim I_{n/2-1}\oplus (C_\phi^{n/2}-1) = I_{n/2-1}\oplus (2^{n/2}-1).\]
Now consider $f_2$. The resultants $\res (f_2,h_1)=2^{n/4}-1$, $\res (f_2,h_2)=1$, $\res (f_2,g_2)=\res (f_2,t^{n/4}-i)\res(f_2,t^{n/4}+i)=2^{n/4}+1$. 
These are pairwise coprime, so Theorem \ref{thm:factorgbetter} implies
\[ f_2(C_g)\sim f_2(C_{h_1})\oplus f_2(C_{h_2})\oplus f_2(C_{g_2}).\]
Moreover, $f_2(C_{h_2})\sim I_{n/4}$. We have $f_2(C_{h_1})=\phi (C_{h_1})$ and $C_\phi=-2$, so Theorem \ref{thm:swaptheorem} implies 
\[\phi(C_{h_1})\sim  I_{n/4-1}\oplus (2^{n/4}-1),\]
and also that $f_2(C_{g_2})\sim I_{n/4}\oplus g_2(C_f)$. In addition,
\[ 
    g_2(C_{f_2}) 
    = C_{f_2}^{n/2}+I_{n/4}
    = (C_{f_2}^{n/4})^2+I_{n/4}
    = (-I_{n/4}-C_{f_2})^2+I_{n/4}
    = \theta (C_{f_2})
\]
where $\theta(t)=2+2t+t^2$. Theorem \ref{thm:swaptheorem} then implies $\theta(C_{f_2})\sim I_{n/4-2}\oplus f_2(C_\theta)$.

Observe that $C_\theta^4=-4I_2$. Hence, 
\[ f_2(C_\theta)= (C_\theta^4)^{n/16} + C_\theta + I = C_\theta + (1 + 2^{n/8}) I \]
 and so, by Theorem \ref{thm:gamma1}, $\gamma_1(f_2(C_\theta))=\cont(t+(1+2^{n/8}))=1$. On the other hand, $\gamma_2(f_2(C_\theta))=\res (f_2,\theta)=\res (f_2,g_2)=2^{n/4}+1$. Thus
   \[ f_2(C_g) \sim I_{n/4}\oplus I_{n/4-1}\oplus (2^{n/4}-1)\oplus I_{n/2-1}\oplus (2^{n/4}+1) \sim I_{n-1}\oplus (2^{n/2}-1).\]
\end{proof}

The cyclically presented groups $G_n(x_0x_1x_{n/2-1})$ ($n$ even), were identified in \cite{EdjvetWilliams,MohamedWilliams} as a particularly challenging subfamily of the family of cyclically presented groups $G_n(x_0x_kx_\ell)$ and were proved in \cite{ChinyereEdjvetWilliams} to be hyperbolic if and only if $n=2,4,6,12,18$. The order and torsion-free rank of their abelianization was calculated in \cite[Theorems 4.1, 4.2]{MohamedWilliams}, and in the finite abelianization case (that is, in the case $\gcd(n,6)=2$), the minimum number of generators of the abelianization was conjectured in \cite[Conjecture 4.4]{MohamedWilliams}. In Theorem \ref{thm:CRSopen} we calculate the structure of the abelianization when finite, proving the conjecture in Corollary \ref{cor:CRSopen}. Theorem \ref{thm:CRSopen} involves the Lucas numbers $L_j$ defined by $L_0=2,L_1=1, L_j=L_{j-1}+L_{j-2}$ ($j\geq 2$), as well as the classical Fibonacci numbers $F_j$ defined in \eqref{eq:FFnumber}. In the proof we make frequent use of the identity (see \cite[p. 200]{NivenZuckermanMontgomery}):
\[ L_n=|\res(t^2+t-1,t^n-1)|+1+(-1)^n.\]

\begin{theorem}\label{thm:CRSopen}
Suppose $\gcd(n,6)=2$, and let $G=G_n(x_0x_1x_{n/2-1})$. Then $G^\mathrm{ab}\cong \Z_{3 L_{n/2}}$, $\Z_{f_{n/4}} \oplus \Z_{15 f_{n/4}}$, $\Z_3\oplus \Z_{L_{n/4}} \oplus \Z_{L_{n/4}}$, $\Z_{L_{n/4}} \oplus \Z_{3 L_{n/4}}$ as $\gcd(n,16)=2,4,8,16$, respectively.
\end{theorem}

For the proof of Theorem \ref{thm:CRSopen} we need the following lemma:

\begin{lemma}\label{lem:Fibonaccigcd}
Suppose $\gcd(n,6)=2$ and $\gcd(n,16)\neq 8$. Then
\begin{itemize}
    \item[(a)] $\gcd(3,L_{n/2}+1+(-1)^{n/2})=1$; and
    \item[(b)] $\gcd(F_{n/2},1+(-1)^{n/2}F_{n/2-1})=1,F_{n/4},L_{n/4}$ if $\gcd(n,16)=2,4,16$, respectively.
\end{itemize}
\end{lemma}

\begin{proof}
\noindent (a) 
Standard results on the Pisano periods of Fibonacci
numbers imply that $3$ divides $(L_{n/2}+1+(-1)^{n/2})$ if and only if $n/2\equiv 4\bmod 8$, or equivalently, $\gcd(n,16)=8$, contrary to hypothesis.

\medskip
 \noindent (b) We calculate $d=\gcd(F_{m},1+\epsilon F_{m-1})$, where $\epsilon=(-1)^{m}$ and $m=n/2$. Similarly to \eqref{eq:fn-j}, an inductive argument shows that for all odd $j$, $1\leq j\leq m/2$,
\begin{alignat}{1}
d=\gcd( F_{m-j}+\epsilon F_j, F_{m-(j+1)}-\epsilon F_{j+1}).
\label{eq:fm-jEpsilonNew}
\end{alignat}
Suppose $\gcd(n,16)=2$. Then $m$ is odd and $\epsilon=-1$. If $j=(m-1)/2$ is odd then substituting this into \eqref{eq:fm-jEpsilonNew} gives \[d=\gcd( F_{(m+1)/2}- F_{(m-1)/2}, F_{(m-1)/2}+ F_{(m+1)/2}=\gcd( F_{(m-3)/2}, F_{(m+3)/2}).\]
If instead $(m-1)/2$ is even then substituting $j=(m+1)/2$ into \eqref{eq:fm-jEpsilonNew} and simplifying yields the same formula:

\[d=\gcd( F_{(m-1)/2}-F_{(m+1)/2}, F_{(m-3)/2}+F_{(m+3)/2})=\gcd( F_{(m-3)/2}, F_{(m+3)/2}).\]
Observe now that, since $\gcd(n,6)=2$ implies $\gcd(m,3)=1$, the subscripts $(m-3)/2$ and $(m+3)/2$, above, are coprime. Using the (standard) property \cite[Theorem II, page 83]{Vajda} $\gcd(F_a,F_b)=F_{\gcd(a,b)}$, we conclude that $d=1$.

Now suppose $\gcd(n,16)=4$. Then $m$ is even, $m/2$ is odd, and $\epsilon=+1$. Substituting $j=m/2$ into \eqref{eq:fm-jEpsilonNew} gives
\[d=\gcd( F_{m/2}+ F_{m/2}, F_{m/2-1}-F_{m/2+1})=\gcd(2F_{m/2},F_{m/2})=F_{m/2}. \]
Finally, suppose $\gcd(n,16)=16$. Then $m$ is even, $m/2$ is even, and $\epsilon=+1$. Substituting $j=m/2-1$ into \eqref{eq:fm-jEpsilonNew} gives
\[d=\gcd( F_{m/2+1}+ F_{m/2-1}, F_{m/2}-F_{m/2})=\gcd( F_{m/2+1}+ F_{m/2-1}, 0)=L_{n/4}.\]
\end{proof}

\begin{proof}[Proof of Theorem \ref{thm:CRSopen}]
Suppose first $\gcd(n,16)\neq 8$. The relation matrix of $G$ is $f(C_g)^T$, where $f(t)=1+t+t^{n/2-1}$, $g(t)=t^n-1=g_1(t)g_2(t)$, where $g_1(t)=t^{n/2}-1, g_2(t)=t^{n/2}+1$. Now
\[\res(f,g_1)=\res(t(1+t+t^{n/2-1}),t^{n/2}-1)=\res(t+t^2+1,t^{n/2}-1)=3\]
(since by hypothesis $\gcd(n/2,3)=1$) so the invariant factors of $f(C_{g_1})$ are 1 ($n/2-1$ times) and 3 (1 time).  Observe further that 
\[\res(f,g_2)= \res(t^2+t-1 , t^{n/2}+1)=L_{n/2}+1+(-1)^{n/2}.\]
By Lemma \ref{lem:Fibonaccigcd} we have $\gcd(3,L_{n/2}+1+(-1)^{n/2})=1$ so, by Theorem \ref{thm:factorgbetter}, the Smith form of $f(C_g)$ is equal to the product of the Smith forms of $f(C_{g_1})$ and $f(C_{g_2})$. Thus it suffices to show that the invariant factors of $f(C_{g_2})$ are $[s_1,s_2]=[1,L_{n/2}],[F_{n/4},5F_{n/4}],[L_{n/4},L_{n/4}]$, as $\gcd(n,16)=2,4,16$, respectively.

Let $h(t)=t^2+t-1$. Then $\res(f,g_2)=\res(h,g_2)$. Moreover, $h(t)\equiv tf(t) \bmod g_2(t)$ so $C_{g_2}f(C_{g_2})=h(C_{g_2})$. Also, $g_2(0)=1$, so $C_{g_2}$ is unimodular, and hence $f(C_{g_2})\sim h(C_{g_2})$. Thus the Smith forms of $f(C_{g_2})$ and of $h(C_{g_2})$ are equal so, in particular, $\gamma_2(g_2(C_{h}))=L_{n/2}+1+(-1)^{n/2}$.

Now $C_{h}=\begin{bmatrix} -1&1\\1& 0 \end{bmatrix}$ and an inductive argument shows that for any $j\geq 1$
\[ C_{h}^j=(-1)^j\begin{bmatrix}F_{j+1} & -F_j\\-F_j& F_{j-1} \end{bmatrix}.\]
Therefore, 
\[ g_2(C_{h})= C_h^{n/2}+I = \begin{bmatrix} 1+(-1)^{n/2}F_{n/2+1} & (-1)^{n/2+1} F_{n/2}\\(-1)^{ n/2+1} F_{n/2}&1+(-1)^{n/2}F_{n/2-1}\end{bmatrix}\]
and hence $\gamma_1(g_2(C_{h}))=d$,  where
\begin{alignat*}{1}
    d
    &=\gcd (1+(-1)^{n/2}F_{n/2+1}, F_{n/2},1+(-1)^{n/2}F_{n/2-1})\\
    &= \gcd(1+(-1)^{n/2}(F_{n/2}+F_{n/2-1}),F_{n/2},1+(-1)^{n/2}F_{n/2-1})\\
    &=\gcd(F_{n/2},1+(-1)^{n/2}F_{n/2-1}).
\end{alignat*}
Then, by Lemma \ref{lem:Fibonaccigcd}, the invariant factors $s_1=\gamma_1,s_2=\gamma_2/\gamma_1$ are given by $[s_1,s_2]=[1,L_{n/2}],[F_{n/4},(L_{n/2}+2)/F_{n/4}],[L_{n/4},(L_{n/2}+2)/L_{n/4}]$, as $\gcd(n,16)=2,4,16$, respectively. If $\gcd(n,16)=4$ then $L_{n/2}+2=5F_{n/4}^2$ \cite[Equation (23)]{Vajda}, and if $\gcd(n,16)=16$ then $L_{n/2}+2=L_{n/4}^2$ \cite[Equation (17a)]{Vajda}. Therefore $[s_1,s_2]=[1,L_{n/2}],[F_{n/4},5F_{n/4}],[L_{n/4},L_{n/4}]$, as $\gcd(n,16)=2,4,16$, respectively, as required.

Now suppose $\gcd(n,16)=8$. Then $n\equiv 8$ or $40\bmod 48$, so $n/4 \equiv 2$ or $6 \bmod 8$ and hence $3 \mid L_{n/4}$. Observe that, working $\bmod~f(t)$, we have $g(t)=(t^{n/2})^2-1\equiv (-(t+t^2))^2-1= h(t)$, where $h(t)=t^4+2t^3+t^2-1=(t^2+t+1)(t^2+t-1)$. Now, applying Theorem \ref{thm:swaptheorem} twice we have
\[ f(C_g)\sim I_{n/2+1}\oplus g(C_f)=I_{n/2+1}\oplus h(C_f)\sim I_{n-4}\oplus f(C_h)\]
so it suffices to obtain the Smith form of $f(C_h)$. By \cite[Lemma 5.1]{NW21}, this is equivalent to computing the Smith form of $f(M)$ where
\[  M=\begin{bmatrix}
   A&0\\
   E&B
\end{bmatrix} \quad \mathrm{where} \quad A=\begin{bmatrix}
    -1&1\\
    1&0
\end{bmatrix}, B=\begin{bmatrix}
    0&-1\\
    1&-1
\end{bmatrix}, E=\begin{bmatrix}
    0&1\\
    0&0
\end{bmatrix}    \]
and hence
\[  f(M)=\begin{bmatrix}
    f(A) & 0\\
    X & f(B)
\end{bmatrix}    \]
where $X$ satisfies
\begin{equation}
XA-BX = Ef(A)-f(B)E.\label{eq:XA-BX}    
\end{equation}
For any $m\geq 1$, $A^m=(-1)^m\begin{bmatrix} F_{m+1}& -F_m\\-F_m&F_{m-1}\end{bmatrix}$, and hence (applying the previous formula for $m=n/2-1$ which is odd) we get $f(A) = \begin{bmatrix}
    -F_{n/2} & 1+F_{n/2-1}\\
    1+F_{n/2-1} & 1-F_{n/2-2}
    \end{bmatrix}$.
We claim that $f(A)= L_{n/4}U$, where 
\[U=\begin{bmatrix}
        -F_{n/4} & F_{n/4-1}\\
        F_{n/4-1} & -F_{n/4-2}
    \end{bmatrix}~\mathrm{is~unimodular~with~inverse}~U^{-1} = \begin{bmatrix}
        F_{n/4-2} & F_{n/4-1}\\
        F_{n/4-1} & F_{n/4}
    \end{bmatrix}.\]
This claim is equivalent to the following Fibonacci and Lucas identities: $F_{n/2}=L_{n/4}F_{n/4}$; $1+F_{n/2-1}=L_{n/4}F_{n/4-1}$; $F_{n/2-2}-1=L_{n/4}F_{n/4-2}$; $F_{n/4}F_{n/4-2}-F_{n/4-1}^2=\epsilon$, where $\epsilon=\pm 1$ (in fact, $\epsilon=(-1)^{n/4-1}$). These can be easily proved using Binet's formula $F_m = (\phi^m - (-\phi)^{-m})/\sqrt{5}$ and the formula $L_m = \phi^m + (-\phi)^{-m}$ where $\phi=(1+\sqrt{5})/2$ is the golden ratio. (Indeed, the first statement is \cite[Equation (13)]{Vajda} and the last statement is Cassini's identity \cite[Equation (29)]{Vajda}). Defining $Y=XU^{-1}$ and noting that $UA=AU$, we thus get
\[ f(M) \sim \begin{bmatrix}
    \gamma I & 0\\
    Y & f(B)
\end{bmatrix} \quad \mathrm{where,~by~\eqref{eq:XA-BX},} \quad   YA-BY = E \gamma - f(B) E U^{-1}. \]
Taking into account that
\[ [(A^T \otimes I) - (I \otimes B) ]^{-1} = \frac12 \begin{bmatrix}
    0&-1&1&0\\
    1&-1&0&1\\
    1&0&1&-1\\
    0&1&1&0
\end{bmatrix},    \]
we can explicitly (and uniquely) solve for $Y$.

We deal with the cases $n \equiv 8,40\bmod 48$ separately. Suppose first $n\equiv 8\bmod 48$. Since $B^3=I$ and $n/2-1\equiv 0\bmod 3$ we have $f(B)=2I+B$, which implies 
\[ f(B)E U^{-1}=  \begin{bmatrix}
        2 F_{n/4-1} & 2 F_{n/4}\\ 
        F_{n/4-1} & F_{n/4}
    \end{bmatrix}.\]
We then obtain
\[  Y = \frac12 \begin{bmatrix}
    3F_{n/4-1}-F_{n/4} & 0\\
    -F_{n/4}-F_{n/4-1} & F_{n/4-1}-F_{n/4}
\end{bmatrix}.   \]
Now
\[ \underbrace{\begin{bmatrix}
    0&1\\
    1&1
\end{bmatrix}}_{:=P} f(B) \begin{bmatrix}
   0&1\\
   1&-1
\end{bmatrix} = \underbrace{1 \oplus 3}_{:=S}, \]
so
\[ f(M) \sim  \begin{bmatrix}
    \gamma I & 0\\
    PY & S
\end{bmatrix}. \] 
Moreover, standard results on the Pisano period of Fibonacci numbers imply that
\[e_2^T PY = -\frac12 \begin{bmatrix}
    2 F_{n/4-2} & F_{n/4-2} 
\end{bmatrix} \equiv 0 \bmod 3,  \]
since $n/4-2\equiv 0\bmod 12$, and thus 
\[ f(M) \sim  \begin{bmatrix}
    \gamma I & 0\\
    0 & S
\end{bmatrix} \sim 1 \oplus 3 \oplus \gamma \oplus \gamma. \] 

Now suppose $n\equiv 40\bmod 48$. In this case $n/2-1\equiv 1\bmod 3$ and so $f(B)=I+2B$, implying 
\[
f(B)EU^{-1} = \begin{bmatrix}
       F_{n/4-1} & F_{n/4}\\
       2 F_{n/4-1} & 2 F_{n/4}       
    \end{bmatrix}. \]
Thus
\[  Y = \frac12 \begin{bmatrix}
    4 F_{n/4-1} & 2 F_{n/4} + F_{n/4-1}\\
    F_{n/4-1} - 2 F_{n/4} & 0
\end{bmatrix}.   \]
Now
\[ \underbrace{\begin{bmatrix}
    0&1\\
    -1&2
\end{bmatrix}}_{:=Q} f(B) \begin{bmatrix}
   0&1\\
   -1&2
\end{bmatrix} = \underbrace{1 \oplus 3}_{:=S}, \]
so
\[ f(M) \sim  \begin{bmatrix}
    \gamma I & 0\\
    QY & S
\end{bmatrix}. \] 
Moreover, again by looking at Pisano periods,
\[e_2^T QY = -\frac12 \begin{bmatrix}
    2 F_{n/4+2} & F_{n/4+2} 
\end{bmatrix} \equiv 0 \bmod 3,  \]
since $n/4+2\equiv 0\bmod 12$, and again we have $f(M)\sim 1 \oplus 3 \oplus \gamma \oplus \gamma$.
\end{proof}

\begin{corollary}[{\cite[Conjecture 4.4]{MohamedWilliams}}]\label{cor:CRSopen}
Suppose $\gcd(n,6)=2$, and let $G=G_n(x_0x_1x_{n/2-1})$. Then
\[d(G^\mathrm{ab}) =\begin{cases}
1 & \mathrm{if}~\gcd(n,16)=2,\\
2 & \mathrm{if}~\gcd(n,16)=4~\mathrm{or}~16,\\
3 & \mathrm{if}~\gcd(n,16)=8.
\end{cases}\]
\end{corollary}

\subsection{Cavicchioli-Repov\v{s}-Spaggiari cyclically presented groups}\label{sec:CRSgroups}

Cavicchioli, Repov\v{s}, Spaggiari  \cite{CRS03} introduced the $8$-parameter family cyclically presented groups
\[ G_n(h,k;m,q;r,s;\ell)=G_n( (x_0x_m\dots x_{m(r-1)})^l(x_hx_{h+q}\dots x_{h+(s-1)q})^{-k})\]
where $n,h,k,m,q,r,s\geq 1$. These form a large class of cyclically presented groups that contain many other well-studied families of cyclically presented groups that arise as fundamental groups of closed, connected, orientable 3-manifolds \cite{CRS03}. The relation matrix of $G_n(h,k;m,q;r,s;\ell)$ is the circulant matrix $f(C_g)^T$ where
\[f(t)= l(t^{rm}-1)/(t^m-1)-kt^h(t^{sq}-1)/(t^q-1),\]
which is of the form $l(t^{rm}-1)/(t^m-1)-kh(t)$ for some $h(t)\in \Z[t]$. In Theorem \ref{thm:lowerboundinvariantfactorsexample} we apply Theorem \ref{thm:genJO} to matrices $f(C_g)$ where $f(t)$ is of this form and $g(t)=t^n-1$. In Corollary \ref{cor:CRSrank} we obtain a lower bound for $d(G^\mathrm{ab})$, the minimum number of generators of the abelianisation of $G=G_n(h,k;m,q;r,s;\ell)$, and record group theoretic and topological consequences.

\begin{theorem}\label{thm:lowerboundinvariantfactorsexample}
Let $f(t)=l(t^{rm}-1)/(t^m-1)-kh(t)$, $g(t)=t^n-1$, where $r,n,m\geq 1, k\in \Z$ with $|k| \neq 1$,  $h(t)\in \Z [t]$. Then $f(C_g)$ has at least 
$\gcd(n,mr)-\gcd(n,m)$ non-unit invariant factors.
\end{theorem}

\begin{proof}
Let $p$ be a prime dividing $k$; at least one such $p$ exists because $|k|\neq 1$. In the notation of Theorem \ref{thm:genJO} we have $f_p(t)=(l \bmod p) \frac{t^{rm}-1}{t^m-1}, g_p(t)=t^n-1 \in \F_p[t]$. If $p$ divides $l$, then $f_p(C_{g_p})=0$ and, by Theorem \ref{thm:genJO}, every invariant factor of $f(C_g)$ is a non-unit, implying the statement since $n > \gcd(n,mr)-\gcd(n,m)$. Assume now $l \not\equiv 0 \bmod p$, and denote $$w_p(t)=\gcd_{\F_p[t]}(f_p(t),g_p(t))\in \F_p[t], \qquad z(t)=\gcd_{\Z[t]}\left(\frac{t^{mr}-1}{t^m-1},g(t)\right)\in \Z[t].$$ Clearly, $z_p(t):=[z(t) \bmod p]\in \F_p[t]$ divides $w_p(t)$, and hence\footnote{We note in passing that, generally, $z_p(t) \neq w_p(t)$. For example, if $n=9$, $m=3$ and $r=p=2$, then $z_2(t)=1$ but $w_2(t)=f_2(t)=t^3+1$.} $\deg w_p(t) \geq \deg z_p(t)=\deg z(t)$.  In turn, by standard results on cyclotomic polynomials,  $\deg z(t)=\sum \varphi(d)$ where $\varphi$ is the Euler totient function and the sum is taken over all integers $d \geq 1$ that divide $n$ and divide $mr$, but do not divide $m$. Hence, by Corollary \ref{cor:genJO2}, the number of non-unit invariant factors of $f(C_g)$ is bounded below by $$ \deg w_p(t) \geq \deg z(t) = \sum_{d|\gcd(n,mr)} \varphi(d) - \sum_{d|\gcd(n,m)} \varphi(d) = \gcd(n,mr) - \gcd(n,m).$$
\end{proof}

\begin{corollary}\label{cor:CRSrank}
Let $n,h,k,m,q,r,s\geq 1$, and let $G=G_n(h,k;m,q;r,s;\ell)$. Then $d(G^\mathrm{ab})\geq \gcd(n,mr)-\gcd(n,m)$. Hence if $G$ is finite then $\gcd(n,mr)-\gcd(n,m)\leq 3$, and if $G$ is solvable then $\gcd(n,mr)-\gcd(n,m)\leq 4$. Moreover, if $G$ is the fundamental group of a closed, connected, orientable 3-manifold $M$ then the Heegaard genus $g(M)\geq \gcd(n,mr)-\gcd(n,m)$.
\end{corollary}

\begin{remark}The lower bound for $d(G^\mathrm{ab})$ in Corollary \ref{cor:CRSrank} is the best possible. To see this consider, for example, the \emph{Sieradski groups} $S(2,n)=G_n(1,1;2,2;2,1;1)$ \cite{Sieradski}. If $12|n$ then $S(2,n)^\mathrm{ab}\cong \Z^2$, so $d(S(2,n)^\mathrm{ab})=2$, which is equal to the lower bound of Corollary \ref{cor:CRSrank}.
\end{remark}

\section*{Acknowledgement}
The second named author was partially supported by the 
Visiting Professor Programme of the Dean of the School of Science at Aalto University, and thanks the Department of Mathematics and Systems Analysis at Aalto University for its hospitality during a research visit in Spring 2024, when this work was completed.

\bibliographystyle{abbrv}
\bibliography{SNF}

\begin{thebibliography}{10}

\bibitem{fastQR}
J.~L. Aurentz, T.~Mach, L.~Robol, R.~Vandebril, and D.~S. Watkins.
\newblock Fast and backward stable computation of roots of polynomials, {P}art
  {II}: {B}ackward error analysis; companion matrix and companion pencil.
\newblock {\em SIAM J. Matrix Anal. Appl.}, 39(3):1245--1269, 2018.

\bibitem{BP94}
D.~A. {Bini} and V.~Y. {Pan}.
\newblock {\em {Polynomial and matrix computations. Fundamental algorithms.
  Vol. 1.}}
\newblock Boston, MA: Birkh\"auser, 1994.

\bibitem{BoileauZieschang}
M.~Boileau and H.~Zieschang.
\newblock Heegaard genus of closed orientable {S}eifert {$3$}-manifolds.
\newblock {\em Invent. Math.}, 76(3):455--468, 1984.

\bibitem{BustomiBarra}
Bustomi and A.~Barra.
\newblock Invertibility of some circulant matrices.
\newblock {\em Journal of Physics: Conference Series}, 893(1):012012, 2017.

\bibitem{CampbellRobertson75}
C.~M. Campbell and E.~F. Robertson.
\newblock {On a class of finitely presented groups of Fibonacci type.}
\newblock {\em {J. Lond. Math. Soc., II. Ser.}}, 11:249--255, 1975.

\bibitem{CRS03}
A.~{Cavicchioli}, D.~{Repov\v{s}}, and F.~{Spaggiari}.
\newblock {Topological properties of cyclically presented groups.}
\newblock {\em {J. Knot Theory Ramifications}}, 12(2):243--268, 2003.

\bibitem{CRS05}
A.~Cavicchioli, D.~Repov\v{s}, and F.~Spaggiari.
\newblock Families of group presentations related to topology.
\newblock {\em J. Algebra}, 286(1):41--56, 2005.

\bibitem{BainsonChinyere}
I.~Chinyere and B.~O. Bainson.
\newblock Perfect {P}rishchepov groups.
\newblock {\em J. Algebra}, 588:515--532, 2021.

\bibitem{ChinyereEdjvetWilliams}
I.~Chinyere, M.~Edjvet, and G.~Williams.
\newblock All hyperbolic cyclically presented groups with positive length three
  relators.
\newblock {\em {P}reprint}, 2024.

\bibitem{ChinyereWilliamsFF}
I.~Chinyere and G.~Williams.
\newblock Fractional {F}ibonacci groups with an odd number of generators.
\newblock {\em Topology and its Applications}, 312:108083, 2022.

\bibitem{Conway65}
J.~H. Conway.
\newblock Advanced problem 5327.
\newblock {\em Amer. Math. Monthly}, 72:915, 1965.

\bibitem{Davis}
P.~J. Davis.
\newblock {\em Circulant matrices}.
\newblock A Wiley-Interscience Publication. John Wiley \& Sons, New
  York-Chichester-Brisbane, 1979.
\newblock Pure and Applied Mathematics.

\bibitem{Horner0}
F.~De~Ter\'an, F.~M. Dopico, and D.~S. Mackey.
\newblock Fiedler companion linearizations and the recovery of minimal indices.
\newblock {\em SIAM J. Matrix Anal. Appl.}, 31(4):2181--2204, 2009/10.

\bibitem{DDPbounds}
F.~De~Ter\'an, F.~M. Dopico, and J.~P\'erez.
\newblock New bounds for roots of polynomials based on {F}iedler companion
  matrices.
\newblock {\em Linear Algebra Appl.}, 451:197--230, 2014.

\bibitem{Dresden}
G.~{Dresden}.
\newblock {Resultants of cyclotomic polynomials.}
\newblock {\em {Rocky Mt. J. Math.}}, 42(5):1461--1469, 2012.

\bibitem{EdjvetWilliams}
M.~Edjvet and G.~Williams.
\newblock The cyclically presented groups with relators {$x_ix_{i+k}x_{i+l}$}.
\newblock {\em Groups Geom. Dyn.}, 4(4):759--775, 2010.

\bibitem{Friedland}
S.~{Friedland}.
\newblock {\em {Matrices. Algebra, analysis and applications.}}
\newblock Hackensack, NJ: World Scientific, 2016.

\bibitem{Hatcher}
A.~{Hatcher}.
\newblock {\em {Algebraic topology.}}
\newblock Cambridge: Cambridge University Press, 2002.

\bibitem{HJ91}
R.~A. Horn and C.~R. Johnson.
\newblock {\em Topics in matrix analysis}.
\newblock Cambridge University Press, Cambridge, 1994.
\newblock Corrected reprint of the 1991 original.

\bibitem{HJ85}
R.~A. Horn and C.~R. Johnson.
\newblock {\em Matrix analysis}.
\newblock Cambridge University Press, Cambridge, second edition, 2013.

\bibitem{Jacobson}
N.~Jacobson.
\newblock {\em Basic algebra. {I}}.
\newblock W. H. Freeman and Company, New York, second edition, 1985.

\bibitem{JohnsonOdoni}
D.~L. {Johnson} and R.~W.~K. {Odoni}.
\newblock {Some results on symmetrically-presented groups.}
\newblock {\em {Proc. Edinb. Math. Soc., II. Ser.}}, 37(2):227--237, 1994.

\bibitem{JohnsonWamsleyWright}
D.~L. {Johnson}, J.~W. {Wamsley}, and D.~{Wright}.
\newblock {The Fibonacci groups.}
\newblock {\em {Proc. Lond. Math. Soc. (3)}}, 29:577--592, 1974.

\bibitem{kap49}
I.~Kaplansky.
\newblock Elementary divisors and modules.
\newblock {\em Trans. Amer. Math. Soc.}, 66:464--491, 1949.

\bibitem{Kim78}
A.~C. Kim.
\newblock Fibonacci varieties.
\newblock {\em Bull. Aust. Math. Soc.}, 19:191--196, 1978.

\bibitem{SteveFiedler}
D.~S. Mackey.
\newblock The continuing influence of {F}iedler's work on companion matrices.
\newblock {\em Linear Algebra Appl.}, 439(4):810--817, 2013.

\bibitem{Maclachlan}
C.~Maclachlan.
\newblock Generalisations of {F}ibonacci numbers, groups and manifolds.
\newblock In {\em Combinatorial and geometric group theory ({E}dinburgh,
  1993)}, volume 204 of {\em London Math. Soc. Lecture Note Ser.}, pages
  233--238. Cambridge Univ. Press, Cambridge, 1995.

\bibitem{MaclachlanReid}
C.~{Maclachlan} and A.~W. {Reid}.
\newblock {Generalised Fibonacci manifolds.}
\newblock {\em {Transform. Groups}}, 2(2):165--182, 1997.

\bibitem{MKS}
W.~Magnus, A.~Karrass, and D.~Solitar.
\newblock {\em Combinatorial group theory}.
\newblock Dover Publications, Inc., Mineola, NY, second edition, 2004.

\bibitem{MohamedWilliams}
E.~Mohamed and G.~Williams.
\newblock An investigation into the cyclically presented groups with length
  three positive relators.
\newblock {\em Experimental Mathematics}, 31(2):537--551, 2022.

\bibitem{NN16}
Y.~{Nakatsukasa} and V.~{Noferini}.
\newblock {On the stability of computing polynomial roots via confederate
  linearizations.}
\newblock {\em {Math. Comput.}}, 85(301):2391--2425, 2016.

\bibitem{NNT17}
Y.~{Nakatsukasa}, V.~{Noferini}, and A.~{Townsend}.
\newblock {Vector spaces of linearizations for matrix polynomials: a bivariate
  polynomial approach.}
\newblock {\em {SIAM J. Matrix Anal. Appl.}}, 38(1):1--29, 2017.

\bibitem{Neuwirth68}
L.~Neuwirth.
\newblock An algorithm for the construction of {$3$}-manifolds from
  {$2$}-complexes.
\newblock {\em Proc. Cambridge Philos. Soc.}, 64:603--613, 1968.

\bibitem{Newmanbook}
M.~Newman.
\newblock {\em Integral matrices}, volume Vol. 45 of {\em Pure and Applied
  Mathematics}.
\newblock Academic Press, New York-London, 1972.

\bibitem{NivenZuckermanMontgomery}
I.~Niven, H.~S. Zuckerman, and H.~L. Montgomery.
\newblock {\em An introduction to the theory of numbers}.
\newblock John Wiley \& Sons, Inc., New York, fifth edition, 1991.

\bibitem{N24}
V.~Noferini.
\newblock Invertible bases and root vectors for analytic matrix-valued
  functions.
\newblock {\em Elec. J. Linear Algebra}, 40:1--13, 2024.

\bibitem{NPer16}
V.~{Noferini} and J.~{P\'erez}.
\newblock {Fiedler-comrade and Fiedler-Chebyshev pencils.}
\newblock {\em {SIAM J. Matrix Anal. Appl.}}, 37(4):1600--1624, 2016.

\bibitem{NP15}
V.~{Noferini} and F.~{Poloni}.
\newblock {Duality of matrix pencils, Wong chains and linearizations.}
\newblock {\em {Linear Algebra Appl.}}, 471:730--767, 2015.

\bibitem{NW21}
V.~Noferini and G.~Williams.
\newblock Matrices in companion rings, {S}mith forms, and the homology of
  3-dimensional {B}rieskorn manifolds.
\newblock {\em J. Algebra}, 587:1--19, 2021.

\bibitem{NW22}
V.~Noferini and G.~Williams.
\newblock Cyclically presented groups as {L}abelled {O}riented {G}raph groups.
\newblock {\em Journal of Algebra}, 605:179--198, 2022.

\bibitem{Horner1}
A.~Ostrowski.
\newblock On two problems in abstract algebra connected with {H}orner's rule.
\newblock In {\em Studies in mathematics and mechanics presented to {R}ichard
  von {M}ises}, pages 40--48. Academic Press, New York, 1954.

\bibitem{Rushanan}
J.~J. Rushanan.
\newblock Eigenvalues and the {S}mith normal form.
\newblock {\em Linear Algebra and its Applications}, 216:177--184, 1995.

\bibitem{SM50}
J.~Sherman and W.~J. Morrison.
\newblock Adjustment of an inverse matrix corresponding to a change in one
  element of a given matrix.
\newblock {\em Ann. Math. Statist.}, 21(1):124--127, 1950.

\bibitem{Sieradski}
A.~J. {Sieradski}.
\newblock {Combinatorial squashings, 3-manifolds, and the third homology of
  groups.}
\newblock {\em {Invent. Math.}}, 84:121--139, 1986.

\bibitem{Smith}
H.~J.~S. Smith.
\newblock {I. On systems of linear indeterminate equations and congruences}.
\newblock {\em Proceedings of the Royal Society of London}, 11:86--89, 1862.

\bibitem{SpaggiariNeuwirth}
F.~Spaggiari.
\newblock A geometric study of generalized {N}euwirth groups.
\newblock {\em Forum Math.}, 18(5):803--827, 2006.

\bibitem{SVNeuwirth}
A.~Szczepa\'{n}ski and A.~Vesnin.
\newblock Generalized {N}euwirth groups and {S}eifert fibered manifolds.
\newblock {\em Algebra Colloq.}, 7(3):295--303, 2000.

\bibitem{Vaidyanathan}
P.~P. Vaidyanathan.
\newblock {\em Multirate Systems and Filter Banks}.
\newblock Prentice Hall, 1993.

\bibitem{Vajda}
S.~Vajda.
\newblock {\em Fibonacci \& {Lucas} numbers, and the golden section. {Theory}
  and applications}.
\newblock Chichester: Ellis Horwood Ltd.; New York: Halsted Press, 1989.

\bibitem{WilliamsCirculantGraphs}
G.~Williams.
\newblock Smith forms for adjacency matrices of circulant graphs.
\newblock {\em Linear Algebra Appl.}, 443:21--33, 2014.

\bibitem{WilliamsLOG}
G.~Williams.
\newblock Generalized {F}ibonacci groups {$H(r,n,s)$} that are connected
  {L}abelled {O}riented {G}raph groups.
\newblock {\em J. Group Theory}, 22(1):23--39, 2019.

\bibitem{Wilson91}
J.~S. Wilson.
\newblock Finite presentations of pro-{$p$} groups and discrete groups.
\newblock {\em Invent. Math.}, 105(1):177--183, 1991.

\bibitem{WilsonNotes}
J.~S. Wilson.
\newblock Finitely presented soluble groups.
\newblock In {\em Geometry and cohomology in group theory ({D}urham, 1994)},
  volume 252 of {\em London Math. Soc. Lecture Note Ser.}, pages 296--316.
  Cambridge Univ. Press, Cambridge, 1998.

\bibitem{Wilson98}
J.~S. Wilson.
\newblock {\em Profinite groups}, volume~19 of {\em London Mathematical Society
  Monographs. New Series}.
\newblock The Clarendon Press, Oxford University Press, New York, 1998.

\end{thebibliography}

\Addresses

\end{document}